\documentclass{elsarticle}

\usepackage{lineno,hyperref}
\modulolinenumbers[5]
\journal{}
\usepackage{ulem}
\usepackage{subfigure}
\usepackage{graphicx}
\usepackage{enumerate}
\usepackage{amsmath,bm}
\usepackage{latexsym}
\usepackage{float}
\usepackage{latexsym}
\usepackage{amsmath}
\usepackage{amssymb}
\usepackage{amsfonts}
\usepackage{verbatim}
\usepackage{mathrsfs}
\usepackage{amsfonts}
\usepackage{colortbl,dcolumn}
\usepackage{amsmath}
\usepackage{amsfonts}
\usepackage{graphicx}
\usepackage{epstopdf}
\usepackage{algorithmic}

\allowdisplaybreaks[4]

\newtheorem{tm}{Theorem}[section]
\newtheorem{rk}{Remark}[section]

\newtheorem{prop}{Proposition}[section]
\newtheorem{lm}{Lemma}[section]
\newtheorem{cor}{Corollary}[section]

\usepackage{enumerate}
\numberwithin{equation}{section}

\newcommand{\E}{\mathbb E}
\newcommand{\PP}{\mathbb P}
\newcommand{\N}{\mathbb N}
\newcommand{\R}{\mathbb R}

\newcommand{\OO}{\mathcal O}

\newcommand{\HH}{\mathbb H}
\newcommand{\LL}{\mathcal L}

\newcommand{\FFF}{\mathscr F}

\newcommand{\<}{\langle}
\renewcommand{\>}{\rangle}
\allowdisplaybreaks[4]

\begin{document}

	\begin{frontmatter}
			\title{Absolute continuity and numerical approximation of stochastic Cahn--Hilliard equation with unbounded noise diffusion}
		\tnotetext[mytitlenote]{This work was supported by National Natural Science Foundation of China (NO. 91530118, NO. 91130003, NO. 11021101, NO. 91630312 and NO. 11290142). } 
		
		\author[cas]{Jianbo Cui\corref{cor}}
		\ead{jcui82@gatech.edu}
		
		\author[cas1]{Jialin Hong}
		\ead{hjl@lsec.cc.ac.cn}
		
		\cortext[cor]{Corresponding author.}
		
		\address[cas]{  School of Mathematics, Georgia Institute of Technology, Atlanta, GA 30332, USA
			}

		\address[cas1]{1. LSEC, ICMSEC, 
			Academy of Mathematics and Systems Science, Chinese Academy of Sciences, Beijing,  100190, China\qquad
			2. School of Mathematical Science, University of Chinese Academy of Sciences, Beijing, 100049, China }

\date{}

\begin{abstract}
In this article, we consider  the absolute continuity and numerical approximation of the solution of the stochastic Cahn--Hilliard equation with unbounded noise diffusion. We first obtain 
the H\"older continuity and Malliavin differentiability of the solution of the stochastic Cahn--Hilliard equation by using the strong convergence of the spectral Gakerkin approximation. Then
we prove the existence and strict positivity of the density function of the law of the exact solution for the stochastic Cahn--Hilliard equation with sublinear growth diffusion coefficient, which fills a gap for the existed result when the diffusion coefficient satisfies a growth condition of order $1/3<\alpha<1$. 
To approximate the density function of the exact solution, 
we 
propose a full discretization based on the spatial spectral Galerkin approximation and the temporal drift implicit Euler scheme. Furthermore, a general framework for deriving the strong convergence rate of the full discretization is developed 
based on the variation approach and the factorization method.  Consequently, we obtain the sharp mean square convergence rates in both time and space via Sobolev interpolation inequalities and semigroup theories.  To the best of our knowledge, this is the first result on the convergence rate of full discretizations for the considered equation. 
\end{abstract}

	\end{frontmatter}

\section{Introduction}

The stochastic  Cahn--Hilliard equation is a fundamental phase field model and can be used to describe the complicated phase separation and coarsening phenomena in a melted alloy that is quenched to a temperature at which only two different concentration phases can exist stably (see e.g. \cite{CH58,NS84}). 
There exist a lot of works focusing on the well-posedness of the stochastic Cahn--Hilliard equation(see e.g. \cite{DD96,Car01,AKM16}).
Recently, the existence and uniqueness of the solution, as well as its regularity estimates, for the stochastic  Cahn--Hilliard equation driven by multiplicative space-time white noise with diffusion coefficient of sublinear growth are proven in \cite{CH19}. However, 
many problems about the solution for the stochastic Cahn--Hilliard equation are still unclear. One of these problems lies on the existence, smoothness, strict positivity of the density of both the exact solution and its approximation. The authors in \cite{Car02} obtain the existence and strict positivity of the density function for stochastic Cahn--Hilliard equation with bounded diffusion coefficient. The authors in \cite{AFK18} prove the absolute continuity of the solution for  stochastic Cahn--Hilliard equation with sublinear growth diffusion coefficient of order $\alpha<\frac 13$. To the best of our knowledge, there still exists no result for the existence of density function for general sublinear growth diffusion coefficient of order $\frac 13\le \alpha<1$. Another important and unclear problem is how to construct the implementary  numerical approximation which can simulate the properties, such as the density function and longtime behaviors, of the exact solution as much as possible. And this topic is far from well-understood since there exists no result on the basic convergent problem of numerical approximations for the the stochastic Cahn--Hilliard equation with unbounded noise diffusion.

In this article, we consider the following stochastic  Cahn--Hilliard equation \begin{align}\label{spde}
dX(t)+A(AX(t)+F(X(t)))dt&=G(X(t))dW(t),\quad t\in (0,T]\\\nonumber
X(0)&=X_0.
\end{align}
Here $0<T<\infty$, $\HH:=L^2(\mathcal O)$ with $\mathcal O=(0,L)$, $L>0,$ 
$-A:D(A)\subset \HH \to \HH$ is the  Laplacian operator under homogenous Dirichlet boundary condition, and  $\{W(t)\}_{t\ge 0}$ is a cylindrical Wiener process on a filtered probability space $(\Omega,\mathcal F,$ 
$\{\mathcal F_t\}_{t\ge 0},\PP)$. 
The nonlinearity $F$ is assumed to be the Nemytskii operator of $f'$, where $f$ is a polynomial of degree 4, i.e., $c_4\xi^4+c_3\xi^3+c_2\xi^2+c_1\xi+c_0$ with $c_i\in \R$, $i=0,\cdots,4$, $c_4>0$.
A typical example is  the double well potential $f=\frac 14(\xi^2-1)^2$. When $G=I$, Eq. \eqref{spde} corresponds to the  Cahn--Hilliard--Cook equation.
The diffusion coefficient $G$ is assumed to be the
Nemytskii operator of $\sigma$, where $\sigma$ is a global Lipschitz function with the sublinear growth condition $|\sigma(\xi)|\le C(1+|\xi|^{\alpha}), \alpha<1$.

To show the existence of the density of the law of the exact solution for Eq. \eqref{spde}, we first consider the spectral Galerkin approximation proposed in \cite{CH19} and derive its H\"older continuity estimates. With the help of truncated  and localization arguments, the Malliavin differentiability of  the spectral Galekin approximation of the truncated equation and the weak Malliavin differentiability of the exact solution are obtained. Combining  with the non-degeneracy assumption on $|\sigma(\cdot)|>0$ and 
the properties of the Green function generated by the  bi-Laplacian operator $-A^2$, we further show that the law of $\{X(t,x_1),X(t,x_2),\cdots,X(t,x_d)\}$ admits a density, with $t>0$, $x_i\in \mathcal O, i=1,\cdots,d$, being different points, which fills the gap on the absolute continuity of the exact solution for stochastic Cahn--Hilliard equation with general sublinear growth diffusion coefficient of order $\frac 13\le \alpha<1$. Furthermore,  assume in addition that $\sigma \in \mathcal C_b^2(\mathbb R)$, by using the transformation on $\Omega$ proposed in \cite{BP98}, the strictly positivity of the density function of the exact solution is shown. 

Once the existence of the density function of the solution is established, an important and natural question is how to simulate this density function numerically. To design such numerical methods to simulate the density function or the distribution of 
the exact solution, it is crucial to solve the basic strong convergent problem. Up to now, there have been a lot of works on the strong convergence problem of numerical approximations for the stochastic Cahn--Hillard equation and it is far from well understood (see e.g. \cite{CCZZ18,FKLL18,KLM11,HJ14,CHS19,QW19}).
The authors in \cite{FKLL18,KLM11}
show the strong convergence of the finite element method and its implicitly full discretization for  Eq. \eqref{spde} driven by additive spatial regular noise. 
The authors in \cite{HJ14} derive the strong convergence rate of the spectral Galerkin method for Eq. \eqref{spde} in one dimension driven by additive  trace class noise.
The convergence rate of the finite element method and its implicitly full discretization is obtained in \cite{QW19} for Eq. \eqref{spde} driven by additive spatial regular noise. 
For Eq. \eqref{spde} driven by additive space-time white noise, we are only aware that the authors in \cite{CHS19} derive the optimal strong convergence rate in space and super-convergence rate in time for a accelerated implicit Euler full  discretization. Recently, the authors in  \cite{CH19} obtain the strong convergence and the optimal strong convergence rate of the the spatial spectral Galerkin method for Eq. \eqref{spde} with unbounded noise diffusion.
However, there exists no result about  the strong convergence of temporally and fully discrete numerical approximations for stochastic Cahn--Hilliard equation with unbounded noise diffusion. 

The present work
considers a full discretization based on the spatial spectral Galerkin method and the temporal drift implicit Euler method, and
makes further contribution on the strong convergence of numerical schemes for stochastic partial differential equations (SPDEs) with non-globally monotone continuous nonlinearity, especially for Eq. \eqref{spde}.  The presence of the unbounded elliptic operator in front of the cubic nonlinearity and the roughness of the driving noise, make the convergence analysis of numerical approximation much more challenging and demanding. 
To overcome such difficulties,  several steps and techniques are introduced.
Instead of studying the strong convergence problem in $\HH$ directly, we use the similar idea in \cite{CHS19} to present the convergence analysis of numerical approximation in 
a negative Sobolev space at the first step. 
To this end, the optimal spatial and temporal regularity estimates of the numerical approximation are presented by transforming Eq. \eqref{spde} to a random PDE and a stochastic convolution, as well as the monotonicity of $-AF$ in $\HH^{-1}$ and the factorization method.
Then the optimal strong convergence error estimate in $\HH^{-1}$ is derived with the help of the variation approach and an appropriate error decomposition.
Based on the Sobolev interpolation inequality and the smoothing effect of the semigroup generated by the bi-Laplacian operator, we recover the optimal convergence rate in mean square sense of the numerical approximation for Eq. \eqref{spde} driven by multiplicative space-time white noise.
Let $\delta t$ be the time stepsize such that $T=K\delta t$, $K\in \N^+$, $  X_0\in \HH^{\gamma}$, $\gamma \in (0,\frac 32)$, $N\in \N^+$ and  $p\ge1$, then
the numerical solution $\{X^N_k\}_{k\le K}$, is strongly convergent to $X$  and  satisfies 
\begin{align*}
\big\|X^N_k-X(t_k)\big\|_{L^2(\Omega;{\HH})}
&\le C(X_0,T,p,\gamma)(\delta t^{(\frac \gamma 2)^-}+\lambda_N^{-(\frac \gamma 2)^-})
\end{align*}
for a positive constant $ C(X_0,T,p,\gamma)$.
We also remark that this approach is also available for deducing the strong convergence rates of numerical schemes for Eq. \eqref{spde}  driven by multiplicative regular noise and for proving the strong convergence of  numerical schemes for Eq. \eqref{spde} under other boundary conditions. 

The outline of this paper is as follows. In the next section,
some preliminaries and assumptions are listed.
Section \ref{sec-3} is devoted to the existence and positivity of the density functions of the laws of the exact solution and its spectral Galerkin approximation, as well as their regularity estimates and Malliavin differentiability.
In Section \ref{sec-pri}, we derive the optimal mean square convergence rate of the implicit full discretization via a new approach based on the variation method and the factorization method.

\section{Preliminaries}\label{sec-2} 
In this section, we give some preliminaries and notations.
Given two separable Hilbert spaces $(\mathcal H, \|\cdot \|_{\mathcal H})$ and $(\widetilde  H,\|\cdot \|_{\widetilde H})$, 
$\LL(\mathcal H, \widetilde H)$ and $\LL_1(\mathcal H, \widetilde H)$  are the Banach spaces of all linear bounded operators 
and  the nuclear operators from $\mathcal H$ to $\widetilde H$, respectively. 
The trace of an operator $\mathcal T\in \LL_1(\mathcal H)$
is $tr[\mathcal T]=\sum_{k\in \N}\<\mathcal Tf_k,f_k\>_{\mathcal H}$, where $\{f_k\}_{k\in \N}$ ($\N=\{0,1,2,\cdots\}$) is any orthonormal basis of $\mathcal H$.
In particular, if $\mathcal T\ge 0$, $tr[\mathcal T]=\|\mathcal T\|_{\mathcal L_1}$.
Denote by $\LL_2(\mathcal H,\widetilde H)$ the space 
of Hilbert--Schmidt operators from $\mathcal H$ into $\widetilde H$, equipped with the usual norm given by  $\|\cdot\|_{\LL_2(\mathcal H,\widetilde H)}=(\sum_{k\in \N}\|\cdot f_k\|^2_{\widetilde H})^{\frac{1}{2}}$.
The following useful property and inequality hold 
\begin{align}\label{tr}
&\|\mathcal S \mathcal T\|_{\mathcal L_2(\mathcal H,\widetilde H)}\le \|\mathcal S\|_{\mathcal L_2(\mathcal H,\widetilde H)}\|\mathcal T\|_{\mathcal L(\mathcal H)},\quad \mathcal T \in \mathcal L(\mathcal H), \;\;\mathcal S\in\mathcal L_2(\mathcal H,\widetilde H),\\\nonumber
&tr[\mathcal Q]=\|\mathcal Q^{\frac 12}\|^2_{L_2(\mathcal H)}=\|\mathcal T\|^2_{ \mathcal L_2(\widetilde H,\mathcal H)},\quad \mathcal Q=\mathcal T \mathcal T^{*},\;\; \mathcal T\in \mathcal L_2(\widetilde H,\mathcal H),
\end{align}
where $\mathcal T^*$ is the adjoint operator of $\mathcal T$.

Given a Banach space $(\mathcal E,\|\cdot\|_{\mathcal E})$, we denote by $\gamma( \mathcal H, \mathcal E)$ the space of $\gamma$-radonifying operators endowed with the norm
$\|\mathcal T\|_{\gamma(\mathcal  H, \mathcal E)}=(\widetilde \E\|\sum_{k\in\N}\gamma_k \mathcal Tf_k \|^2_{\mathcal E})^{\frac 12}$,
where $(\gamma_k)_{k\in\N}$ is a Rademacher sequence on a
probability space $(\widetilde \Omega,\widetilde \FFF, \widetilde \PP)$.
For convenience, let  $L^q:=L^q(\OO)$, $2\le q<\infty$ and 
 equipped with the usual inner product and norm.
We also need the following BurkerH\"older inequality in $L^q$,
\begin{align}\label{Burk}
\left\|\sup_{t\in [0,T]}\Big\|\int_0^t \phi(r)dW(r)\Big\|_{L^q}\right\|_{L^p(\Omega)}
&\le 
C_{p,q}\|\phi\|_{L^p(\Omega;  L^2([0,T]; \gamma(\mathcal H;L^q))}\\\nonumber 
&\le  C_{p,q}\Big(\E\Big(\int_0^T\Big\|\sum_{k\in \N} (\phi(t) f_k)^2\Big\|_{L^{\frac q2}}dt\Big)^{\frac p2}\Big)^{\frac 1p}, 
 \end{align}
 where $\{f_k\}_{k\in \N}$ is any orthonormal basis of $\mathcal H$.

Next, we introduce some assumptions and spaces associated with $ A$.
 We denote by 
$H^k:=H^k(\mathcal O)$ the standard Sobolev space and $E:=\mathcal C(\overline{\mathcal O})$ the continuous function space. 
Denote $A=-\Delta $ the Dirichlet Laplacian operator with 
$$D(A)=\left\{v\in H^2(\mathcal O): v =0\;\; \text{on} \;\; \partial \mathcal O\right\}.$$
It is known that $A$ is a positive definite, self-adjoint and unbounded linear operator on $\HH$. Thus there exists an orthonormal 
eigensystem $\{(\lambda_j,e_j)\}_{j\in \N}$ such that $0<\lambda_1\le \cdots \le \lambda_j\le \cdots$ with $\lambda_j \sim j^{ 2}.$
We  define 
$\HH^{\alpha}$, $\alpha\in \R$ as the space of the series $v:=\sum_{j=1}^{\infty}v_je_j$, $v_j\in \R$, such that
$\|v\|_{\HH^{\alpha}}:=(\sum_{j=1}^{\infty}\lambda_j^{\alpha}v_j^2)^{\frac 12}<\infty$. Equipped with the norm $\|\cdot\|_{\HH^{\alpha}}$ and corresponding inner product, the Hilbert space $\HH^{\alpha}$ equals $D(A^{\frac \alpha 2})$. Throughout this article, we denote by $C$ a generic constant which may depend on several parameters but never on the projection parameter $N$ and may change from occurrence to occurrence.

To study the density function of the solution at fixed $t\in[0,T], x\in \mathcal O$, some objects of Malliavin calculus are also introduced. Denote by $\mathcal S$ as the 
space of simple functionals, i.e., $\mathbb F=f(W(h_1), \cdots, W(h_m))$ with $f\in \mathbb C^{\infty}_p(R^m)$ and an orthonormal sequence $\{h_i\}_{i=1}^m$ in $\mathscr H:=L^2(0,T;\mathbb H).$ Here $$W(h):=\int_{0}^T h(s)dW(s)$$ for $h\in \mathscr H$.
For $\mathbb F\in \mathcal S$, the Malliavin derivative of $\mathbb F$ is $\mathcal H$-valued random variable
$$D_{t,x}\mathbb F=\sum_{i=1}^m \frac {\partial f}{\partial x^i}(W(h_1),W(h_2),\cdots,W(h_m))h_i(t,x).$$
Similarly, the $k$th derivative of $\mathbb F$ is given by 
$$D_{\alpha}^k\mathbb F=\sum_{i_1,\cdots, i_k=1}^m\partial_{i_1}\cdots\partial_{i_k}f(W(h_1),\cdots,W(h_m))h_{i_1}(\alpha_1) \cdots h_{i_k}(\alpha_k),$$
where $\alpha:=(\alpha_1,\cdots,\alpha_k),$ $\alpha_i=(t_i,x_i)\in [0,T]\times \mathcal O$.

For $p\ge 1$ and $k\in \mathbb N^+$, the space $\mathcal D^{p,k}$
is the completion of $\mathcal S$ with respect to the seminorm $$\|\mathbb F\|_{\mathcal D^{p,k}}=\|\mathbb F\|_{L^{p}(\Omega;\mathbb R)}+\sum_{i=1}^k\| D^i \mathbb F\|_{L^{p}(\Omega; \mathscr H^{\otimes k})}.$$ We also define $\mathcal D^{\infty}=\cap_{p\ge 1}\cap_{k\in \mathbb N^+}\mathcal D^{p,k}.$ The covariance matrix associated to a $\mathbb R^d$-valued random vector $\mathbb F$ is denoted by $C_{ij}=\<D^1\mathbb F^i,D^1\mathbb F^j\>$, for $ i,j=1,\cdots,d.$

\section{Absolute continuity of the exact solution}
\label{sec-3}

In this section, we first give the well-posedness result of the consider equation \eqref{spde}, as well as the regularity estimate and strong convergence of the spectral Galerkin approximation. 
Next, we study the Malliavin differentiability of the exact solution via the spectral Galerkin approximation of the truncated equation of \eqref{spde}.  
Then under the non-degeneracy condition of the diffusion coefficient, we prove that the absolute continuity of the exact solution with respect to the Lebsgue measure. We would like to mention that all the results on the density function of this section still hold if the boundary condition of Laplacian operator is changed to the homogenous Neumann condition.

\subsection{Well-posedness and Galerkin approximation of the stochastic Cahn--Hilliard equation}

In this subsection, we first present the well-posedness result of the considered equation  whose proof can be found in \cite{CH19}.

\begin{tm}\label{lm-well}
Let $T>0$, $ X_0\in \HH^{\gamma}, \gamma\in(0,\frac 32)$, $p\ge1$. Then
Eq. \eqref{spde} possesses a unique mild solution
$X$ in $L^{p}(\Omega; C([0,T]; \HH))$.
Moreover,   there exists $ C(X_0,T,p,\gamma)>0$  such that 
\begin{align}\label{reg-x}
\big\|X\big\|_{L^p(\Omega;C([0,T];\HH^{\gamma}))}
&\le C(X_0,T,p,\gamma)
\end{align}
and 
\begin{align}\label{reg-tm}
\|X(t)-X(s)\|_{L^p(\Omega;\HH)}
&\le C(X_0,T,p,\gamma)(t-s)^{\frac \gamma  4},
\end{align}
where $s,t\in[0,T]$.
\end{tm}

\begin{rk}\label{con-co}
Under the condition of Theorem \ref{lm-well} with $\gamma\ge \frac 32$, by using the Sobolev embedding theorem and the factorization method, we can get that for $\beta<1,$
\begin{align*}
\|X \|_{L^p(\Omega;\mathcal C([0,T];\mathcal C^{\beta}(\mathcal O)))}<\infty
\end{align*}
and for $\beta<\frac 38$,
\begin{align*}
\|X \|_{L^p(\Omega;\mathcal C^{\beta}([0,T]; E))}<\infty.
\end{align*}
\end{rk}

Denote $P^N$  the spectral Galerkin projection into the linear space  spanned by the first $N$ eigenvectors  $\{e_1,\cdots,e_N\}$.
Then the spectral Galerkin approximation satisfies the following SPDE
\begin{align}\label{sga}
dX^N(t)+A(AX^N(t)+P^NF(X(t)))dt&=P^NG(X^N(t))dW(t),\quad t\in (0,T],\\\nonumber
X^N(0)&=P^NX_0.
\end{align}

Recall that the following strong convergence property in \cite{CH19} of the spectral Galerkin approximation.

\begin{lm}\label{lm-xn}
Let $X_0\in \HH^{\gamma}$, $\gamma\in (0,\frac 32)$ and $p\ge1$.
There exists a unique solution  $X^N$  of Eq. \eqref{sga}  satisfying 
\begin{align}\label{p-mom}
\big\|X^N\big\|_{L^p(\Omega;C([0,T];\HH^{\gamma}))}
&\le C(X_0,T,p,\gamma)
\end{align}
and 
\begin{align*}
\|X^N(t)-X^N(s)\|_{L^p(\Omega;\HH)}
&\le C(X_0,T,p,\gamma)(t-s)^{\frac \gamma  4},
\end{align*}
where $C(X_0,T,p,\gamma)>0$ and $s,t\in [0,T]$.

If in addition assume that $\sup\limits_{N\in \N^+}\|X^N_0\|_{E}\le C(X_0)$, then for $\beta<\gamma$, there 
exists $ C(X_0,T,p,\beta)>0$ such that 
\begin{align}\label{strong-spa}
\big\|X^N-X\big\|_{L^p(\Omega;C([0,T];{\HH}))}
&\le C(X_0,T,p,\beta)\lambda_N^{-\frac \beta 2}.
\end{align}
\end{lm}

In the following, we first show that the exact solution of \eqref{spde} is absolutely continuous with respect to the Lebesgue measure if in addition assume that
$\sigma$ is non-degenerate, i.e., 
\begin{align*}
|\sigma(x)|\ge c_0>0, \; \text{for} \; x\in \mathbb R
\end{align*}
and 
that $\sigma$ is continuously differentiable on $\mathbb R$.
Then to study the strict positivity of the joint distribution of $(X(t,x_1), X(t,x_2),\cdots, X(t,x_d))$, we further assume that $\sigma\in \mathcal C^{2}$ and possesses bounded and Lipschitz derivatives.  

For simplicity, throughout this rest part of section 3, we assume that $X_0$ is a sufficient smooth function. We also remark that if 
$|\sigma(x)|>0$, then by the similar arguments in \cite{BP98}, all these results in this section still hold. 

\subsection{Malliavin differentiability of the solution}

Our strategy to obtain the existence of the density function of $X(t,x)$ is split into three steps:
first, we truncated the \eqref{spde} by adding a cut-off function $\theta_R$ on the nonlinear drift part and study the properties of truncated equation by using its semi-discretization. Then we show that the Malliavin differentiability of $X_R(x,t)$ via the spectral Galerkin approximation, which implies the local
Malliavin differentiability of $X(x,t).$ Finally, 
by using a localization argument and a lower bound estimate of the Malliavin matrix, we obtain  the  absolute continuity of the exact solution.

Let $\theta_R: \mathbb R \to \mathbb R^+$ be a $\mathcal C^{\infty}$ cut-off function satisfying $\theta_R\le 1$ and possessing bounded derivatives.  
For any $ R$,
$\theta_R(x):=1,$ if $|x|\le  R$, and $\theta_R(x):=0,$ if $|x|> R+1$. Then we consider the truncated version of \eqref{spde},
\begin{align}\label{trun-spde}
dX_R&=-A^2X_Rdt-A(F_R(X_R))dt+G(X_R)dW(t),\\
X_R&=X_0, 
\end{align}
where $F_R$ is the Nemytskii operator of the function $\theta_R(\cdot)f'(\cdot).$ 

By denoting an increasing sample subset sequence $\{\Omega_R\}_{R\in \N^+}$,$$\Omega_R=\{\omega\in \Omega: \sup\limits_{t\in [0,T]}\|X(t,x)\|_E\le R\},$$
due to the existence of continuous modification of  $X$ by Remark \ref{con-co}, 
we have that $\lim\limits_{R\to \infty} P(\Omega_R)= P(\Omega)=1$
and that 
$(X_R,\Omega_R)$ is a localization of $(X,\Omega)$.

\begin{lm}\label{pri-xr}
There exists a unique solution $X_R$ of \eqref{trun-spde} 
which satisfies
\begin{align*}
\E[\sup_{t\in [0,T]}\|X_R\|^{p}_{E}]\le C(R,X_0,T,p), 
\end{align*}
for $p\ge 2$.
\end{lm}

\textbf{Proof.}
We use the spectral Galerkin method to semi-discretize \eqref{trun-spde} and define its solution as $X_R^N$.
Since $F_R$ is Lipschitz, by the same arguments in \cite{CH19}, we can get that for $p\ge2$, $\E[\sup\limits_{t\in [0,T]}\|X_R^N\|^{p}_{E}]$ and that 
$$\lim\limits_{N\to\infty}\E[\sup_{t\in [0,T]}\|X_R^N-X_R\|^{p}_{E}] =0,$$
which completes the proof.
\qed

Next, we show that $X_R(t,x) \in \mathcal D^{1,2}$ for $R\in \mathbb R^+$, $t\in [0,T], x\in \mathcal O$, which implies that $X(t,x)\in \mathcal D_{loc}^{1,2}$. For convenience, 
we denote $G(t,x,y)$ as the Green function determined by the double Laplacian operator under homogeneous Dirichlet condition. In this case, $G(t,x,y)=\sum_{i=1}^{\infty}e^{-\lambda_i^2t}e_i(x)e_i(y),$ for $x,y\in \mathcal O$. Denote $G^N(t,x,y)=\sum_{i=1}^{N}e^{-\lambda_i^2t}e_i(x)e_i(y),$ for $x,y\in \mathcal O$ and $N\in \mathbb N^+$.

\begin{prop}\label{mal-xr}
For fixed $(t,x)\in [0,T]\times\mathcal O$, $X_R(t,x)\in \mathcal D^{1,2}$ and satisfies 
\begin{align}\label{mal-xr}
D_{s,y}X_R(t,x)&=G(t-s,x,y)\sigma(X_R(s,y))\\\nonumber
&\quad+\int_s^t\int_{\mathcal O}\Delta G(t-r,x,z)\Big(F_R'(X_R(r,z)) D_{s,y}X_R(r,z)\Big)dzdr\\\nonumber
&\quad+\int_s^t\int_{\mathcal O}G(t-r,x,z)\sigma'(X_R(r,z))D_{s,y}X_R(r,z)W(dz,dr),
\end{align}
where $s\in [0,t).$
\end{prop}

\textbf{Proof.}
Consider the spectral Galerkin approximation of $X_R^N$ which satisfies the following stochastic differential equation
\begin{align*}
dX^N_R=-A^2X^N_Rdt-A P^NF_R(X^N_R)dt+P^N\sigma(X^N_R)dW(t).
\end{align*}
Then the Gateaux derivative $D_{s}^{e_i}X^N_R$ for a basis $\{e_i\}_{i\in \N^+}$ of $L^2(\mathcal O)$ exists and satisfies
\begin{align*}
dD_{s}^{e_i}X^N_R&=-A^2D_{s}^{e_i}X^N_Rdt-A P^NF_R'(X^N_R)D_{s}^{e_i}X^N_Rdt
+P^N\sigma'(X^N_R)D_{s}^{e_i}X^N_RdW(t),\\
D_{s}^{e_i}X^N_R(s)&=P^N\sigma(X^N_R)e_i.
\end{align*}
Then we can get that $D_{s,y}X_R^N(t,x)=\sum_{i\in \N^+}\<D_{s,\cdot}X_R^N(t,x),e_i\>e_i(y)$ formally satisfies 
\begin{align}\label{mal-xrn}
D_{s,y}X_R^N(t,x)&=G^N(t-s,x,y)\sigma(X_R^N(s,y))\\\nonumber 
&\quad+\int_s^t\int_{\mathcal O}\Delta G(t-r,x,z)P^N\Big(F_R'(X_R(r,z)) D_{s,y}X_R^N(r,z)\Big)dzdr\\\nonumber 
&\quad+\int_s^t\int_{\mathcal O}G(t-r,x,z)P^N\sigma'(X_R(r,z))D_{s,y}X_R^N(r,z)W(dz,dr).
\end{align}
To show that $D_{s,y}X_R^N(t,x)$ is the solution of the above integration equation, we need to show that the map defined by the right hand side of the above equality has a unique fixed point. By a rescaling argument, it sufficient to give a priori estimate on $D_{s,y}X_R^N(t,x)$.\\

Taking $L^2$-norm with respect to $y$, we get 
\begin{align*}
&\| D_{s,y}X_R^N(t,x)\|_{L^2_y}^2\\
&\le 3 \|G^N(t-s,x,y)\sigma(X_R(s,y))\|_{L^2_y}^2\\
&\quad+3\|\int_s^t\int_{\mathcal O}\Delta G(t-r,x,z)P^N\Big(F_R'(X_R(r,z)) D_{s,y} X_R(r,z) \Big)dzdr\|_{L^2_y}^2\\
&\quad+3\|\int_s^t\int_{\mathcal O}G(t-r,x,z)P^N\sigma'(X_R(r,z))D_{s,y} X_R(r,z)W(dz,dr)\|_{L^2_y}^2\\
&:=IV_1+IV_2+IV_3.
\end{align*}
By taking $L^1([0,T])$ integration and the property that $D_{s,y} X_R^N(r,z)=0$ if $s>r$,
we get 
\begin{align*}
\int_0^T\|D_{s,y}X_R^N(t,x)\|_{ L^2_y}^2ds&=\int_0^t\|D_{s,y}X_R^N(t,x)\|_{ L^2_y}^2ds\\
&\le \int_0^tIV_1ds+\int_0^t IV_2ds+\int_0^tIV_3ds
\end{align*}
Due to Lemma \ref{pri-xr} and the sublinear growth of $\sigma$, we have that 
\begin{align*}
\E\big[\int_0^TIV_1ds\big]&\le C\int_0^T\|G^N(t-s,x,y)\sigma(X_R^N(s,y))\|_{L^2_y}^2ds\\
&\le C(R)<\infty.
\end{align*}
The H\"older inequality, the boundedness of the truncated  drift term $F_R$, the smoothy estimate of $G(t,x,y)$ and the Fubini's theorem yield that 
\begin{align*}
\E\big[\int_0^t IV_2ds\big]&\le 
C\E \Big[\int_0^t\|\int_s^t\int_{\mathcal O}\Delta G(t-r,x,z)P^N\Big(F_R'(X_R^N(r,z)) \\
&\qquad D_{s,y}X_R^N(r,z)\Big)dzdr\|_{L^2_y}^2ds\Big]\\
&\le C\int_0^t \int_{\mathcal O}\E \Big[\Big(\int_s^t|\int_{\mathcal O}\Delta G^N(t-r,x,z))\\
&\qquad F_R'(X_R^N(r,z))D_{s,y} X_R^N(r,z)dz|dr\Big)^2 \Big]dyds
\\
&\le \int_0^t \int_{\mathcal O}\E \Big[\Big(\int_s^t\int_{\mathcal O} |\Delta G^N (t-r,x,z))|_{L_x^{\infty}} |F_R'(X_R^N(r,z))|^2dzdr\\
&\quad \times
\int_s^t\int_{\mathcal O}  |\Delta G^N(t-r,x,z))|_{L_x^{\infty}} 
|D_{s,y} X_R^N(r,z)|^2dzdr\Big) \Big]dyds\\
&\le C(R) \int_0^t \int_s^t  |\Delta G^N (t-r,x,z))|_{L_x^{\infty}} 
 \int_{\mathcal O}\E\Big[\|D_{s,y} X_R^N(r,z)\|_{L^2_y}^2\Big] dzdrds\\
 &\le C(R) \int_0^t (t-r)^{-\frac 12} \int_0^r 
 \int_{\mathcal O}\E\Big[\|D_{s,y} X_R^N(r,z)\|_{L^2_y}^2\Big] dzdsdr
\end{align*} 

The Burkh\"older inequality, Fubini's theorem, the boundedness of $\sigma'$ and the form of $G^N$ imply that 
\begin{align*}
&\E[\int_0^tIV_3ds]\\
&\le 
C \int_0^t\E\Big[\|\int_s^t\int_{\mathcal O}G^N(t-r,x,z)\sigma'(X_R(r,z))D_{s,y} X_R(r,z)W(dz,dr)\|_{L^2_y}^2\Big]ds\\
&\le C\int_0^t \int_{\mathcal O}\E\Big[\int_{s}^t (t-s)^{-\frac 1 2 }\|\sigma'(X_R(r,z))\|_{L^{\infty}_z}^2\int_{\mathcal O}\|D_{s,y} X_R(r,z)\|^2 dzdr\Big]dyds\\
&\le C\int_0^t (t-r)^{-\frac 12} \int_{0}^r\int_{\mathcal O}\E\Big[\|D_{s,y} X_R(r,z)\|_{L^2_y}^2\Big] dzds dr.
\end{align*}
Thus by the Gronwall inequality, we conclude that  
$\sup\limits_{N}\sup\limits_{t\in [0,T]}\sup\limits_{x\in \mathcal O}\|X_R^N(t,x)\|_{\mathcal D^{1,2}}\le C(R)$ and $X_R^N$ satisfies \eqref{mal-xrn}.
By the strong convergence of $X^N_R$, we have that $X_R(t,x)\in \mathcal D^{1,2},$
 and that $D_{s,y}X_R(t,x)$ satisfies \eqref{mal-xr}.
\qed
 
 In fact, due to the strong convergence of $X^N_R$, we can further show that $X^N_R$ is strongly convergent to 
 $X_R$ in $\mathcal D^{1,2}$ for any fixed $t>0,x\in \mathcal O$. 
 Based on the above proposition, we obtain the following weak differentiability of the exact solution $X\in \mathcal D_{loc}^{1,2}$. 
 The Malliavin derivative $D_{s,y}X$ is defined by the restriction on $\Omega_R$.  
 \begin{cor}
 For any $t>0, x\in \mathcal O$.
 The solution $X(t,x)$ of \eqref{spde} belongs to $\mathcal D_{loc}^{1,2}$.
 \end{cor}

 Similar arguments yield that $X^N_R(t,x)\in \mathcal D^{1,p}, p\ge 2$
 and $X_R(t,x)\in \mathcal D^{1,p}, p\ge 2$ and that  $X(t,x)\in \mathcal D_{loc}^{1,p}$.  With a slight modification, one can obtain that $$\|D_{s,y}X_R(t,x)\|_{L^p(\Omega;L^2_y)}\le C(R,T,X_0)(t-s)^{-\frac 14}.$$

\subsection{Absolute continuity of the law of the exact solution}

Assume that $x_i\in \mathcal O, i=1,\cdots,d$ are distinct points.
In this part, we show that $\mathbb F_R=(X_R(t,x_1), X_R(t,x_2),\cdots, X_R(t,x_d))$ possesses a density function on $\mathbb R^d$.
Without loss of generality, we assume that $x_i$ is strictly increasing with respect to $i$.

To this end, according to Lemma 2.3.1 in \cite{Nua06}, it suffices  to show that the Malliavin 
matrix $C^R$ of $\mathbb F_R$, defined by $C^R_{ij}(t)=\<D\mathbb F_R^i(t),D\mathbb F_R^j(t)\>, t>0$, $i,j\le d$, is invertible and satisfies that for $p\ge 2$, there exists a small parameter $\epsilon_0(p)$ such that for all $\epsilon\le \epsilon_0(p)$,
\begin{align*}
\sup_{|\xi|=1}\mathbb P(\xi^TC^R(t)\xi\le \epsilon)\le \epsilon^p. 
\end{align*}

\begin{prop}\label{abo-xr}
For $p\ge 2$ and $T>0$, there exists $\epsilon_0(p)$ such that for all $\epsilon\le \epsilon_0(p)$, $t\in [0,T],$
\begin{align*}
\sup_{|\xi|=1}\mathbb P(\xi^TC^R(t)\xi\le \epsilon)\le \epsilon^p. 
\end{align*}
\end{prop}

\textbf{Proof.}
Denote $D_{s,y}X_R(t,x)=G(t-s,x,y)\sigma(X_R(s,y))+Q_{s,y}(t,x)$, where 
\begin{align*}
Q_{s,y}(t,x)X_R(t,x)&=\int_s^t\int_{\mathcal O}\Delta G(t-r,x,z)F'_R(X_R(r,z))
D_{s,y}X_R(r,z)drdz\\
&\quad+
\int_s^t \int_{\mathcal O}G(t-r,x,z)\sigma'_R(X_R(r,z))
D_{s,y}X_R(r,z)W(dz,dr).
\end{align*}
We use the method similar to \cite{BP98} and show the lower bound $\xi^TC^R(t)\xi$, with $|\xi|=1$. Assume that $\min_{1\le i\le d-1}|x_i-x_{i+1}|\ge L>4\epsilon^\alpha$ and that there exists a parameter $\alpha\in (\frac 14, 1)$ such that $(x_i-\epsilon^{\alpha},x_i+\epsilon^{\alpha})\subset \mathcal O.$
Then we get 
\begin{align*}
\xi^TC^R(t)\xi &\ge \sum_{i=1}^d\int_{t-\epsilon}^{t}\int_{x_i-\epsilon^{\alpha}}^{x_i+\epsilon^{\alpha}}
(\sum_{j=1}^dD_{r,z}X_R(t,x_j)\xi_j)^2dzdr\\
&\ge 
\frac 12 \sum_{i=1}^d\int_{t-\epsilon}^{t}\int_{x_i-\epsilon^{\alpha}}^{x_i+\epsilon^{\alpha}}
D_{r,z}X_R(t,x_i)^2\xi_i^2dzdr\\
&\quad-\sum_{i=1}^d\int_{t-\epsilon}^{t}\int_{x_i-\epsilon^{\alpha}}^{x_i+\epsilon^{\alpha}}
(\sum_{j\neq i}^dD_{r,z}X_R(t,x_j)\xi_j)^2dzdr\\
&=:\frac 12V_1(\epsilon)-V_2(\epsilon).
\end{align*}

We first give an upper estimate of $V_2(\epsilon)$.
Since $\{(x_i-\epsilon^{\alpha},x_i+\epsilon^{\alpha})\}_{i=1}^d$ are disjointed interval, it holds that 
\begin{align*}
\sup_{|\xi|=1}V_2(\epsilon) &\le C_d\sum_{j\neq i}^d\int_{t-\epsilon}^{t}\int_{x_i-\epsilon^{\alpha}}^{x_i+\epsilon^{\alpha}}
(\sigma(X_R(r,z))G(t-r,x_j,z))^2dzdr\\
&\quad +C_d\sum_{j\neq i}^d\int_{t-\epsilon}^{t}\int_{x_i-\epsilon^{\alpha}}^{x_i+\epsilon^{\alpha}}
(Q_{r,z}(t,x_j))^2dzdr\\
&=: C_d\sum_{j\neq i}^d V_{21}^{ij}(\epsilon)+C_d\sum_{j\neq i}^d V_{22}^{ij}(\epsilon).
\end{align*}

By taking $p$-moment on both sides of the above inequality, and using H\"older's inequality, the sublinear growth of $\sigma$,  and 
the property of the Green function that 
\begin{align*}
|G(t,x,y)|\le C_2t^{-\frac 14} \exp(-C_1\frac {|x-y|^{\frac 43}}{t^{\frac 13}}),
\end{align*}
we get for $i\neq j$,
\begin{align*}
\E\Big[|V_{21}^{ij}(\epsilon)|^p\Big]
&\le C_{p,d} \E\Big[\Big|\int_{t-\epsilon}^{t}\int_{x_i-\epsilon^{\alpha}}^{x_i+\epsilon^{\alpha}}
\sigma^2(X_R(r,z))(t-r)^{-\frac 12}e^{-2C_1 L^{\frac 43}\epsilon^{-\frac 13}}dzdr\Big|^p\Big]\\
&\le C_{p,d}\epsilon^{p(\alpha +\frac 12)}e^{-C\epsilon^{-\frac {1}{3}}}.
\end{align*}
By the definition of $Q_{r,z}(t,x)$, we obtain 
\begin{align*}
&\E\Big[|V_{22}^{ij}(\epsilon)|^p\Big]\\
&\le C_{p,d}\E\Big[\Big|\int_{t-\epsilon}^{t}\int_{x_i-\epsilon^{\alpha}}^{x_i+\epsilon^{\alpha}}
(Q_{r,z}(t,x_j))^2dzdr\Big|^p\Big]\\
&\le C_{p,d} \E\Big[\Big(\int_{t-\epsilon}^{t}\int_{x_i-\epsilon^{\alpha}}^{x_i+\epsilon^{\alpha}}
\Big|\int_r^t\int_{\mathcal O}\Delta G(t-\theta,x_j,z_1)\\
&\quad F'_R(X_R(\theta,z_1))
D_{r,z}X_R(\theta,z_1)d\theta dz_1\Big|^2drdz\Big)^p\Big]\\
&\quad+
 C_{p,d}\E\Big[\Big(\int_{t-\epsilon}^{t}\int_{x_i-\epsilon^{\alpha}}^{x_i+\epsilon^{\alpha}}\Big|\int_r^t\int_{\mathcal O}G(t-r,x_j,z_1)\sigma'_R(X_R(\theta,z_1))\\
 &\quad 
D_{r,z}X_R(\theta,z_1)W(dz_1,d\theta)\Big|^2drdz\Big)^p\Big]\\
&=: V_{221}^{ij}(\epsilon)+ V_{222}^{ij}(\epsilon).
\end{align*}
H\"older's inequality and the estimate of the Green function yield that 
\begin{align*}
V_{221}^{ij}(\epsilon)
&\le C\epsilon^{(\alpha+1)(p-1)}
\int_{t-\epsilon}^{t}\int_{x_i-\epsilon^{\alpha}}^{x_i+\epsilon^{\alpha}}\E\Big[\Big(\int_r^t\int_{\mathcal O}|\Delta G(t-\theta,x_j,z_1)|\\
&\quad |F'_R(X_R(\theta,z_1))|
|D_{r,z}X_R(\theta,z_1)|d\theta dz_1\Big)^{2p}\Big]drdz,\\
&\le 
C\epsilon^{(\alpha+1)(p-1)}
\int_{t-\epsilon}^{t}\int_{x_i-\epsilon^{\alpha}}^{x_i+\epsilon^{\alpha}}
\Big(\int_r^t (t-\theta)^{-\frac 12}\\
&\quad 
\Big(\sup_{z_1}\E\big[\|D_{r,z}X_R(\theta,z_1)\|_{L^2_z}^{2p}\big]\Big)^{\frac 1{2p}} d\theta \Big)^{2p}drdz\\
&\le C \epsilon^{(\alpha+1)p+\frac p2}.
\end{align*}
Similar arguments, together with the Burkh\"older inequality,
imply that 
\begin{align*}
V_{222}^{ij}(\epsilon)
&\le C\epsilon^{(\alpha+1)(p-1)}\int_{t-\epsilon}^{t}\int_{x_i-\epsilon^{\alpha}}^{x_i+\epsilon^{\alpha}}
\E\Big[\Big|\int_r^t\int_{\mathcal O}G(t-r,x_j,z_1)\sigma'_R(X_R(\theta,z_1))\\
 &\quad 
D_{r,z}X_R(\theta,z_1)W(dz_1,d\theta)\Big|^{2p}\Big]dzdr\\
&\le C \epsilon^{(\alpha+1)p+\frac p4}.
\end{align*} 
Thus we conclude that $\E[|\sup_{|\xi|=1}V_2(\epsilon)|^p]\le C \epsilon^{p(\alpha +1)+\frac p4}$.

Next, we estimate the lower bound of the term $\sup_{|\xi|=1}V_1(\epsilon)$.
By Young inequality, we have that 
\begin{align*}
\sup_{|\xi|=1}V_1(\epsilon)
&\ge 
\frac 12 \sum_{i=1}^d\int_{t-\epsilon}^{t}\int_{x_i-\epsilon^{\alpha}}^{x_i+\epsilon^{\alpha}}
\sigma^2(X_R(r,z))G^2(t-r,x_i,z)\xi_i^2dzdr\\
&\quad-
\sum_{i=1}^d\int_{t-\epsilon}^{t}\int_{x_i-\epsilon^{\alpha}}^{x_i+\epsilon^{\alpha}}
Q^2_{r,z}(t,x_i)\xi_i^2dzdr\\
&=:\frac 12 \sum_{i=1}^dV_{11}^{ii}(\epsilon) - \sum_{i=1}^d V_{12}^{ii}(\epsilon).
\end{align*}
By repeating the similar procedures of estimating $V_{22}^{ij}(\epsilon)$, we get that 
\begin{align*}
\E\big[|V_{12}^{ii}(\epsilon)|^p\big]&\le \epsilon^{(\alpha+1)p+\frac p4}. 
\end{align*}
For the term $V_{11}^{ii}(\epsilon)$, we have 
\begin{align*}
&V_{11}^{ii}(\epsilon)\\
&\ge \frac 12\int_{t-\epsilon}^{t}\int_{x_i-\epsilon^{\alpha}}^{x_i+\epsilon^{\alpha}}
\sigma^2(X_R(r,z))G^2(t-r,x_i,x_i)\xi_i^2dzdr\\
&\quad-\int_{t-\epsilon}^{t}\int_{x_i-\epsilon^{\alpha}}^{x_i+\epsilon^{\alpha}}
\big(\sigma(X_R(r,z))G(t-r,x_i,x_i)-\sigma(X_R(r,z))G(t-r,x_i,z)\big)^2\xi_i^2dzdr\\
&=:\frac 12V_{111}^{ii}(\epsilon)- V_{112}^{ii}(\epsilon).
\end{align*}

Based on the lower bound estimate of the Green function 
\begin{align*}
G(s,x_j,x_j)\ge C_{x_j}s^{-\frac 14},
\end{align*}
and the fact that $x_i \in \mathcal O$ implies that $\min_{i=1}^dC_{x_i}>0$,  together with the non-degeneracy of $\sigma$, we have that 
\begin{align*}
V_{111}^{ii}(\epsilon)&\ge C\inf_{x}{|\sigma(x)|}\int_{t-\epsilon}^{t}\int_{x_i-\epsilon^{\alpha}}^{x_i+\epsilon^{\alpha}}
(t-r)^{-\frac 12}dzdr \xi_i^2\\
&\ge C\epsilon^{\frac 12+\alpha}\xi_i^2.
\end{align*}
Thus we conclude that 
$\sum_{i=1}^dV_{111}^{ii}\ge C\epsilon^{\frac 12+\alpha}.$

Then we give an upper bound estimate of $V_{112}^{ii}$ by using the continuity of $X_R$. 
The H\"older inequality and the estimate of Green function 
\begin{align*}
|G(s,x,y)-G(s,x,y_1)|\le C|y-y_1|s^{-\frac 12}
\end{align*}
yield that for $p\ge 2$,
\begin{align*}
&\E \Big[|V_{112}^{ii}(\epsilon)|^p\Big]\\
&\le C\E\Big[\Big|\int_{t_0-\epsilon}^{t_0}\int_{x_i-\epsilon^{\alpha}}^{x_i+\epsilon^{\alpha}}
\big|\sigma(X_R(r,z))\big|^2\big(G(t-r,z,x_i)-G(t-r,x_i,x_i)\big)^2dzdr\Big|^p\Big]\\
&\le C\E\Big[\Big|\int_{t_0-\epsilon}^{t_0}\int_{x_i-\epsilon^{\alpha}}^{x_i+\epsilon^{\alpha}}
\big|\sigma(X_R(r,z))\big|^2|z-x_i|(t-r)^{-\frac 34}dzdr\Big|^p\Big]\\
&\le C\epsilon^{(\frac 14+2\alpha)p}.
\end{align*}
Combining the above estimates, we get that 
\begin{align*}
\xi^TC^R(t)\xi
\ge 
\frac 12V_1(\epsilon)-V_2(\epsilon)
&\ge
C_{p,d,R}\epsilon^{(\frac 12+\alpha)}-
\Big(V_2(\epsilon)+\sum_{i=1}^d(V_{12}^{ii}(\epsilon)+V_{112}^{ii}(\epsilon))\Big).
\end{align*}
Thus the Chebyshev inequality yields that \begin{align*}
&\sup_{|\xi|=1}\mathbb P(\xi^TC^R(t)\xi\le \epsilon)\\
&\le
\sup_{|\xi|=1}\mathbb P\Big(V_2(\epsilon)+\sum_{i=1}^d(V_{12}^{ii}(\epsilon)+V_{112}^{ii}(\epsilon))
\le \epsilon+C_{p,d,R}\epsilon^{\frac 12+\alpha}\Big)\\
&\le (\epsilon+C_{p,d,R}\epsilon^{\frac 12+\alpha})^{-p}\E\Big[\Big|V_2(\epsilon)+\sum_{i=1}^d(V_{12}^{ii}(\epsilon)+V_{112}^{ii}(\epsilon))\Big|^p\Big]\\
&\le (\epsilon+C_{p,d,R}\epsilon^{\frac 12+\alpha})^{-p} C_{p,d,R} \epsilon^{(2\alpha+\frac 14)p},
\end{align*}
which completes the proof.
\qed

\begin{tm}\label{den-x}
Let $t>0$, $x_i\in \mathcal O, i=1,\cdots,d$ be different points. Then the law of $\{X(t,x_1),X(t,x_2),\cdots,X(t,x_d)\}$ admits a density.
\end{tm}

\textbf{Proof.}
Based on Proposition \ref{abo-xr} and the localization arguments, we have that for any $p\ge 2$, there exists $R_0$  large enough, such that  for $R\ge R_0$,
\begin{align*}
\sup_{|\xi|=1}\mathbb P\Big(\xi^TC(t)\xi=0\Big)
&=\sup_{|\xi|=1}\mathbb P\Big(\xi^TC(t)\xi=0, \omega\in \Omega_{R}^c\Big)\\
&\quad+\sup_{|\xi|=1}\mathbb P\Big(\xi^TC^R(t)\xi =0,\omega\in \Omega_R\Big)\\
&\le \epsilon+\sup_{|\xi|=1}\mathbb P\Big(\xi^TC^R(t)\xi=0\Big)
\\
&\le \epsilon.
\end{align*}
Let $R\to \infty$, we get 
\begin{align*}
\sup_{|\xi|=1}\mathbb P\Big(\xi^TC(t)\xi=0\Big)=0,
\end{align*}
which, together with $X(t,x)\in \mathcal D_{loc}^{1,2}$, implies 
the existence of the density function of the joint distribution.
\qed

Indeed, the above argument also yields the existence and smoothness of the density function of the numerical approximation. 
By using the localization $X^N_K:=X^NI_{\Omega_K},$ 
$$\Omega_K=\{\omega\in \Omega: \sup_{t\in[0,T]}\sup_{x\in\mathcal O} |X^N(t,x)|\le K\},$$ and using similar arguments in proving Theorem \ref{den-x}, we can get the following result.

\begin{cor}
Let $t>0$, $x_i\in \mathcal O, i=1,\cdots,d,$ be different points. Then the law of $\{X^N(t,x_1),X^N(t,x_2),\cdots,X^N(t,x_d)\}$ admits a density.
\end{cor}

\subsection{Positivity of the density function}

To illustrate the strategy to show the positivity of the density function, we use the idea of \cite{BP98} and denote the probability space $(\Omega,\mathcal F, \mathbb P)$ as $\Omega=\{\phi\in C([0,T]\times \mathcal O)|\phi(0,x)=\phi(t,0)=0\}$ such that the canonical process is a Brownian sheet. 
Consider a sequence $\{h_{\epsilon_n}\in L^2([0,T]\times \mathcal O)\}$, $\epsilon_n\to 0$ and $z\in \mathbb R^d.$ Let $T_n^z$ be the transformation on $\Omega$ such that 
\begin{align*}
T_n^z(\omega)(t,x)
=\omega(t,x)+\sum_{i=1}^dz_i\int_0^t\int_{\mathcal O}h_{\epsilon_n}^i(y,s)dyds.
\end{align*}
Define $\Phi_n(z):=\mathbb F\circ T_n^z$.
The function $h_{\epsilon_n}$ is chosen according to the lower bound estimate in Proposition \ref{abo-xr}, where $\epsilon_n$  satisfies the same condition of $\epsilon$ in Proposition \ref{abo-xr} for $n$ large enough.  Let $\alpha\in (\frac 14,1)$. 
Then we define 
\begin{align*}
h_{\epsilon_n}^{i}(s,y)=(c_n^{i})^{-1}1_{[t-\epsilon_n,t]}(s)1_{(x_i-\epsilon_n^{\alpha},x_i+\epsilon_n^{\alpha})}(y),
\end{align*}
where 
\begin{align*}
c_n^{(i)}=\int_{t-\epsilon_n}^t\int_{x_i-\epsilon_n^{\alpha}}^{x_i+\epsilon_n^{\alpha}}G(t-r,x_i,y)dyds.
\end{align*}
Then it suffices to prove that for some fixed $R$, there will exist $C,C'$ and $\delta$, $h_{\epsilon_n},$  $r$ small enough,
\begin{align*}
&\mathbb P\Big(\Omega_R \cap \Big\{|\mathbb F-y|\le r \Big\}\cap \Big\{\det(\nabla \Phi_n(0))\ge C\Big\}
\\
&\qquad \cap \Big\{\sup_{|z|\le \delta}\Big(|\nabla \Phi_n(z)|+|\nabla^2 \Phi_n(z)|\Big)\le C'\Big\}\Big)>0,
\end{align*}
where  $\Phi_n(z)=\mathbb F\circ T_n^z.$

The above result will be immediately obtained if the following two propositions hold.

\begin{prop}\label{prop-est-1}
Let $y\in supp(\mathbb P\circ \mathbb F^{-1})$.
There exists $C,r_0>0$ such that for $r\le r_0$,
\begin{align*}
\limsup_{n\to \infty}\mathbb P\Big(\Omega_R \cap \Big\{|\mathbb F-y|\le r \Big\}\cap \Big\{\det(\nabla \Phi_n(0))\ge C\Big\}\Big)>0.
\end{align*}
\end{prop}

\textbf{Proof.}
For simplicity, we only give the proof for the case $d=1$.
Denote  the Malliavin derivative of $X_R(s,x)\circ T_z^n$ with respect to the direction $h_n$ as
$u_z^R(s,x)=\frac {\partial}{\partial z} X_R(s,x)\circ T_z^n.$
Then we have 
\begin{align*}
u_{z}^R(t,x)&=\int_0^t\int_{\mathcal O}G(t-s,x,y)\sigma(X_R(s,y)\circ T_z^n)h_{\epsilon_n}(s,y)dyds\\
&\quad +\int_{t-\epsilon}^t\int_{\mathcal O}\Delta G(t-s,x,y)F_R'(X_R(s,y)\circ T_z^n)u_{z}^R(s,y)dyds\\
&\quad +\int_{t-\epsilon}^t\int_{\mathcal O}G(t-s,x,y)\sigma'(X_R(s,y)\circ T_z^n)u_{z}^R(s,y)W(dy,ds)\\
&\quad +\int_{t-\epsilon}^t\int_{\mathcal O}G(t-s,x,y)\sigma'(X_R(s,y)\circ T_z^n)u_{z}^R(s,y)zh_{\epsilon_n}(s,y)dyds.
\end{align*}
For the first term in the above equality, 
\begin{align*}
&\E\Big[\Big|\int_0^t\int_{\mathcal O}G(t-s,x,y)\sigma(X_R(s,y)\circ T_z^n)h_{\epsilon_n}(s,y)dyds\Big|^p\Big]\\
&\le C_{p,R}\sup_{s\in [0,t]}\E\Big[|\sigma(X_R(s,y)\circ T_z^n)|^p_{L^{\infty}_y}\Big]\le C_{p,R}<\infty,
\end{align*}
where we use $\sup_{s\in[0,T]}\E\Big[|X_R(s,y)\circ T_z^n|_{L^{\infty}_y}^p\Big]\le C_p$ which can be proven by similar procedures in the estimate of $X_R(s,y).$
By the similar arguments in Proposition \ref{mal-xr}, we get that for $p$ large enough, $$\E [ |u_{z}^R(t,x)|^p]\le C_{p,R}<\infty,$$
and  that 
\begin{align*}
&\E [ \Big|u_{z}^R(t,x)-\int_0^t\int_{\mathcal O}G(t-s,x,y)\sigma(X_R(s,y)\circ T_z^n)h_{\epsilon_n}(s,y)dyds\Big|^p\Big]\\
&\le C_{p,R}(\epsilon_n^{\frac {p}4}+|z|^p).
\end{align*}
Notice that 
\begin{align*}
&\int_0^t\int_{\mathcal O}G(t-s,x,y)\sigma(X_R(s,y)\circ T_0^n)h_{\epsilon_n}(s,y)dyds\\
&=\sigma(X_R(t,x))+
\int_0^t\int_{\mathcal O}G(t-s,x,y)\Big(\sigma(X_R(s,y)-\sigma(X_R(t,x))\Big)h_{\epsilon_n}(s,y)dyds.
\end{align*}
Taking $z=0$ and using the continuity of $X_R(\cdot,\cdot)$, we obtain that  
\begin{align*}
\E\Big [ \Big|u_{0}^R(t,x)-\sigma(X_R(t,x))\Big|^p\Big]
&\le C_{p,R}\epsilon_n^{\frac {p}4}.
\end{align*}
Thus for any $y \in supp(\mathbb P\circ \mathbb F^{-1})$, there exists $r_0>0$ such that for $0<r\le r_0$,
$\mathbb P\Big(|\mathbb F-y|\le r \Big)>0.$
Furthermore, for $r\le r_0$, it holds that
$$\limsup_{n\to \infty}\mathbb P\Big(\Omega_R\cap\Big\{|\mathbb F-y|\le r\Big\} \cap  \Big\{\det(\nabla_z \Phi_n(0)) \ge C\Big\}\Big)>0.$$
\qed

To prove the following proposition, we in addition assume that $|\sigma''|_{L^{\infty}(\mathbb R)} \le C$.

\begin{prop}\label{prop-est-2}
There exists $C'>0$ and small enough $\delta >0$ such that 
\begin{align*}
\lim_{n\to \infty} \mathbb P\Big( \Big\{\sup_{|z|\le \delta}\Big(|\nabla \Phi_n(z)|+|\nabla^2 \Phi_n(z)|\Big)\ge C'\Big\}=\mathcal O(\delta).
\end{align*}

\end{prop}

\textbf{Proof.}
For simplicity, we only show the estimate of $|\nabla \Phi_n(z)|$ for  $d=1$. For the estimate of $|\nabla \Phi_n(z)|$ with $d>1$ and the estimate of $|\nabla^2 \Phi_n(z)|$, the proofs are similar.
Denote 
\begin{align*}
v_z^R(t,x)&=\int_0^t\int_{\mathcal O}G(t-s,x,y)\sigma(X_R(s,y)\circ T_z^n)h_{\epsilon_n}(s,y)dyds\\
&\qquad+\int_{t-\epsilon}^t\int_{\mathcal O}G(t-s,x,y)\sigma'(X_R(s,y)\circ T_z^n)u_{z}^R(s,y)zh_{\epsilon_n}(s,y)dyds.
\end{align*}
Based on the estimates in the proof of Proposition \ref{prop-est-1}, taking $\delta$ small enough, we have that 
\begin{align*}
\E\Big[|u_z^R(t,x)-u_{z_1}^R(t,x)|^p\Big]+\E\Big[|v_z^R(t,x)-v_{z_1}^R(t,x)|^p\Big]&\le C_{p, R}|z-z_1|^p.
\end{align*}
By using the Kolmogorov criterion of continuity theorem, we obtain
the H\"older continuity of $u_z^R(t,x)$ and $v_z^R(t,x)$,
and that 
\begin{align*}
&\E\Big[\sup_{|z|\le \delta }|u_z^R(t,x)|^p\Big]+\E\Big[\sup_{|z|\le \delta }|v_z^R(t,x)|^p\Big]\le C_{p,R}.
\end{align*}
Similarly, we have that 
\begin{align*}
&\lim_{n\to \infty}\sup_{|z|\le \delta }\E\Big [ \Big|u_{z}^R(t,x)-v_z^R(t,x)\Big|^p\Big]=0,
\end{align*}
which, together with the Garsia-Rodemich-Rumsey lemma, yields that  
\begin{align*}
\lim_{n\to \infty}\E\Big [\sup_{|z|\le \delta } \Big|u_{z}^R(t,x)-v_z^R(t,x)\Big|^p\Big]=0.
\end{align*}
Now, it suffices to show that  
\begin{align}\label{supz-pro=0}
&\lim_{n\to \infty} \mathbb P\Big( \Big\{\sup_{|z|\le \delta}\Big|\int_0^t\int_{\mathcal O}G(t-s,x,y)\\\nonumber
&\qquad \qquad \qquad \sigma(X_R(s,y)\circ T_z^n)h_{\epsilon_n}(s,y)dyds\Big| \ge C'\Big\}\cap \Omega_R \Big)=0.
\end{align}
Since on $\Omega_R$, 
 $$\sup_{s\in[0,T]}\sup_{y\in\mathcal O}|\sigma(X_R(s,y))|\le C(1+|X_R(s,y)|^{\alpha})\le C(1+R^{\alpha}),\alpha<1,$$
we have that
\begin{align*}
&\Big|\int_0^t\int_{\mathcal O}G(t-s,x,y)\sigma(X_R(s,y)\circ T_z^n)h_{\epsilon_n}(s,y)dyds\Big|\\
&\le \Big|\int_0^t\int_{\mathcal O}G(t-s,x,y)\sigma(X_R(s,y))h_{\epsilon_n}(s,y)dyds\Big|\\
&\quad+\Big|\int_0^t\int_{\mathcal O}G(t-s,x,y)(\sigma(X_R(s,y)\circ T_z^n)-\sigma(X_R(s,y))h_{\epsilon_n}(s,y)dyds\Big|\\
&\le C(1+R)+\Big|\int_0^t\int_{\mathcal O}G(t-s,x,y)(\sigma(X_R(s,y)\circ T_z^n)-\sigma(X_R(s,y))h_{\epsilon_n}(s,y)dyds\Big|.
\end{align*}
By the H\"older continuity of $X_R(s,y)\circ T_z^n$
with respect to $z$ and using Markov inequality, we can choose proper $C'(R)$ such that 
\eqref{supz-pro=0} holds.
\qed

Based on the above propositions, we obtain the following theorem, which  yields that $supp(\mathbb P\circ \mathbb F^{-1})=\mathbb R^d$, see e.g. \cite{BP98}.

\begin{tm}
For $t>0$, distinct points $x_i\in \mathcal O$, $i=1,\cdots,d$, there exists a continuous function 
$\rho>0$ such that for all $\phi\in \mathcal C_b^{\infty}(\mathbb R^d,\mathbb R^+)$,
\begin{align*}
\E\Big[\phi(X(t,x_1),\cdots,X(t,x_d))\Big]
\ge \int_{\mathbb R^d}\phi(y)\rho(y)dy.
\end{align*} 
\end{tm}
\textbf{Proof.}
Based on Propositions \ref{prop-est-1} and 
\ref{prop-est-2}, by the similar proof of Theorem 3.3 in \cite{BP98},
we have
\begin{align*}
\E\Big[1_{\Omega_R}\phi(X_R(t,x_1),\cdots,X_R(t,x_d))\Big]
\ge \int_{\mathbb R^d}\phi(y)\rho_R(y)dy,
\end{align*}
which, together with $\phi(\cdot)\in \mathbb R^+$, completes the proof.
\qed

\section{Regularity estimate and strong convergence of a full discretization}
\label{sec-pri}

In previous section, we establish the existence and positivity of the density function of the exact solution  $X$. In the viewpoint of applications, it is important to design numerical approximation for computing the density function of $X$. However,  it  is far from well-understood since there exists no result on the basic convergent problem of numerical approximations for the stochastic Cahn--Hilliard equation with unbounded noise diffusion. 
In this section, we aim to solve the  basic strong convergence problem in terms of  the following full discretization.

Let $\delta t$ be the time stepsize such that $T=K\delta t$, $K\in \N^+$.  
The full discrete numerical scheme starting from $X^N_0:=P^NX_0$, based on the temporal drift implicit Euler scheme  and the spatial spectral Galerkin method, is define as 
\begin{align}\label{full0}
X^N_{k+1}&=X^N_k-A^2X^N_{k+1}\delta t
-AP^NF(X^N_{k+1})+G(X_{k}^N) \delta_k  W, \; k\le K-1
\end{align}  where $ \delta_k  W=W(t_{k+1})-W(t_k), \; k\le K-1$.
Then we also have the mild form of $X^N_k$,
\begin{align}\label{full}
X^N_{k+1}=T_{\delta t} X^N_{k}
-\delta tT_{\delta t}A
P^N(F(X^N_{k+1}))+T_{\delta t}G(X_{k}^N) \delta_k  W, \; k\le K-1,
\end{align}
where $T_{\delta t}=(I+A^2\delta t)^{-1}$.

\subsection{A priori and regularity estimates of the full discretization}

In order to study the moment estimate of the numerical scheme, we decompose  $X^N_k$ into $Y^N_k$ and $Z^N_k$, which satisfy for $k\le K-1$ and $N\in \N^+$,
\begin{align}\label{ily}
Y^{N}_{k+1}&=Y^N_{k}-A^2Y^N_{k+1}\delta t-AF(Y^N_k+Z^N_k)\delta t, \\\label{ilz}
Z^N_{k+1}&=Z^N_{k}-A^2Z^N_{k+1}\delta t+G(Y^N_k+Z^N_k) \delta_k  W,
\end{align}
with $Y^N_0=X^N_0, Z^N_0=0.$
We would like to mention that the similar decomposition has been used to derive the strong convergence rates of numerical approximations for non-global Lipschitz SPDEs driven by additive noise(see e.g. \cite{BJ16,BCH18,CH18}). 
The following lemmas and propositions are
 devoted to deducing the a priori estimates of $Y^{N}_k$, $Z^{N}_k$ and $X^N_k$, $k\le K$, $N\in \N^+$.

\begin{lm}
For $2\le p\le \infty$ and  $m > 0$, it holds that 
\begin{align*}
\|T_{\delta t}^m f\|_{L^p}\le C(m\delta t)^{-\frac 14(\frac 12 -\frac 1p)}\|f\|, \quad f\in \HH,\; t>0.
\end{align*}
\end{lm}

\textbf{Proof.}
For $f\in \HH$, we have $T_{\delta t}^mf=\sum\limits_{i\in \N^+}(\frac 1{1+\lambda_i^2\delta t})^{m}\<f,e_i\>e_i$.
From the uniform boundedness of $e_i$, $i\in \N^+$ and the fact that $ci^2 \le \lambda _i\le Ci^{2}$, it follows that 
\begin{align*}
\|T_{\delta t}^mP^Nf\|_{E}&=\|\sum_{i=1}^{N}(\frac 1{1+\lambda_i^2\delta t})^{m}\<f,e_i\>e_i\|_{ E}\\
& \le (\sum_{i=1}^{N}(\frac 1{1+\lambda_i^2\delta t})^{2m})^{\frac 12}\|f\|\le C(m\delta t)^{-\frac 18}\|f\|,
\end{align*}
and 
\begin{align*}
\|T_{\delta t}^mP^Nf\|&\le \|\sum_{i=1}^{N}(\frac 1{1+\lambda_i^2\delta t})^{m}\<f,e_i\>e_i\|\le C\|f\|.
\end{align*}
The Riesz--Thorin interpolation theorem leads to the desired result. 
\qed
\begin{lm}\label{pri-xnk}
Let $X_0\in \HH$, $T>0$ and $q\ge1$.
There exists a unique solution  $\{X^N_k\}_{k\le K}$  of Eq. \eqref{full0}  satisfying 
\begin{align}\label{dis-p-mom}
\sup_{k\le K}\E\Big[ \big\|X^N_k\big\|_{\HH^{-1}}^q\Big]
&\le C(X_0,T,q)
\end{align}
for a positive constant $ C(X_0,T,q)$.
\end{lm}

\textbf{Proof.}
Taking inner product with $ Y^N_{k+1}$ in $\HH^{-1}$ on both sides of Eq. \eqref{ily} leads to
\begin{align*}
\|Y^N_{k+1}\|_{\HH^{-1}}^2
\le \|Y^N_{k}\|_{\HH^{-1}}^2
-2\|\nabla  Y^N_{k+1}\|^2\delta t
-2\<F(Y^N_{k+1}+Z^N_{k+1}), Y^N_{k+1}\>\delta t.
\end{align*}
From the monotonicity of $-F$, the equivalence of norms in $\HH^1$ and $H^1$, and the Young inequality, it follows that for some small $\epsilon>0$,
\begin{align}\label{pri-h-11}
&\| Y^N_{k+1}\|_{\HH^{-1}}^2
+8(c_4-\epsilon)\| Y^N_{k+1}\|_{L^4}^4\delta t+
(2-\epsilon)\| Y^N_{k+1}\|_{\HH^{1}}^2\delta t\\\nonumber 
&\le 
\| Y^N_{k}\|_{\HH^{-1}}^2
+C(\epsilon)\Big(1
+\|Z_{k+1}^N\|_{L^{4}}^{4}+\|A^{-\frac 12} Z_{k+1}^N\|^2\Big)\delta t\\\nonumber 
&\le 
\| Y^N_{k}\|_{\HH^{-1}}^2
+C(\epsilon)\Big(1
+\|Z_{k+1}^N\|_{L^{4}}^{4}\Big)\delta t.
\end{align}
From the mild form of $Z^N_k$, the H\"older inequality  and the Burkh\"older inequality, it follows that for $p\ge 2$ and $q\ge 4$, 
\begin{align*}
&\E[\|Z^N_{k+1}\|_{L^{p}}^q]\\
&=
\E\Big[\big\|\sum_{i=0}^{k} T_{\delta t}^{k+1-i}G(Y^N_{i}+Z^N_i) \delta  W_i\big\|_{L^p}^{q}\Big]\\
&\le 
\E\Big[\Big(\sum_{i=0}^{k}\sum_{j\in \N^+}\big\| T_{\delta t}^{k+1-i}G(Y^N_{i}+Z^N_i)e_j\big\|_{L^p}^2\delta t\Big)^{\frac q2}\Big]\\
&\le \E\Big[\Big(\sum_{i=0}^{k}
((k+1-i)\delta t)^{-\frac 12+\frac 1{2p}}(1+\|Y^N_{i}\|^{2\alpha}+\|Z^N_i\|^{2\alpha})\delta t\Big)^{\frac q2}\Big]\\
&\le C(\sum_{i=0}^{k}
((k+1-i)\delta t)^{-1+\frac 1{p}})^{\frac q4} \E\Big[\Big(\sum_{i=0}^{k}(1+\|Y^N_{i}\|^{4\alpha}+\|Z^N_i\|^{4\alpha})\delta t\Big)^{\frac q4}\Big].
\end{align*}
Using the Young inequality, we obtain for $0\le s\le t$,
\begin{align*}
\E[\|Z^N_{k+1}\|_{L^{p}}^q]
&\le  CT^{\frac q{4p}} \Big(1+\E[(\sum_{i=0}^{k}\|Y^N_i\|^{4 \alpha}\delta t)^{\frac q4}]+\sum_{i=0}^k\E[\|Z_{i}^N\|^{\alpha q}]\delta t\Big)
\\ 
&\le CT^{\frac q{4p}} \Big(1+\E[(\sum_{i=0}^{k}\|Y^N_i\|^{4 \alpha}\delta t)^{\frac q4}]+\sum_{i=0}^k\E[\|Z_{i}^N\|^{q}]\delta t\Big).
\end{align*}
Gronwall's inequality yields that for  $0\le s\le t\le T$,
\begin{align}\label{pri-znk-lp}
\E[\|Z^N_{k+1}\|_{L^{p}}^q]
&\le C(T)\Big(1+\E[(\sum_{i=0}^{k}\|Y^N_i\|^{4 \alpha}\delta t)^{\frac q4}]\Big).
\end{align}

Now taking $m$-moment, $m\in \N^+$ on \eqref{pri-h-11} and letting $p=4$,$q=4m,$ we have 
\begin{align*}
&\E\Big[(\sum_{i=0}^k\|Y^N_k\|^4_{L^4}\delta t)^m\Big]\\
&\le C\|Y^N(0)\|_{\HH^{-1}}^{2m}
+C(\epsilon,T)\sum_{i=0}^k(1+\E[\|Z^N_i\|_{L^4}^{4m}])\delta t\\
&\le C\|Y^N(0)\|_{\HH^{-1}}^{2m}+C(\epsilon,T)\sum_{i=0}^k\Big(C(\epsilon_1)+\epsilon_1 \E[(\sum_{j=0}^i \|Y^N_k\|_{L^4}^{4}\delta t)^{m}]\delta t\Big),
\end{align*}
where $\epsilon_1>0$ is a small number such that 
$C(\epsilon,T)\epsilon_1T<\frac 12$. 
The above estimation leads to 
\begin{align*}
\E\Big[(\sum_{i=0}^k\|Y^N_k\|^4_{L^4}ds)^m\Big]
&\le  C\|Y^N(0)\|_{\HH^{-1}}^{2m}+C(m,\epsilon,\epsilon_1,T),
\end{align*}
which in turn yields that for $m\in \N^+$,
\begin{align}\label{dis-pri-h-1}
&\E\Big[\|Y^N_{k+1}\|^{2m}_{\HH^{-1}}\Big]+\E\Big[(\sum_{i=0}^k\|\nabla Y^N_k\|^2\delta t)^m\Big]+\E\Big[(\sum_{i=0}^k\|Y^N_k\|^4_{L^4}\delta t)^m\Big]\\\nonumber
&\le C(X_0,T,m),
\end{align}
and for $p\ge 2$, $q\ge 4$,
\begin{align}\label{pri-znk0}
\E\Big[\|Z^N_{k+1}\|^{q}_{L^p}\Big]
&\le C(X_0,T,m).
\end{align}
Combining the above estimates, we complete the proof.
\qed
\begin{lm}\label{sup-e-znk}
Let  $X_0\in \HH$, $T>0$ and $q\ge1$.
There exists a positive constant $ C(T,q,X_0)$
such that 
\begin{align*}
\sup_{k\le K}\E\Big[\|Z^N_k\|^{q}_{E}\Big]
&\le C(T,q,X_0).
\end{align*}
\end{lm}
\textbf{Proof.}
By the Burkholder inequality and the smoothing effect of $T_{\delta t}$, we have that for $ q \ge 2$,
\begin{align*}
\E\Big[\|Z_{k+1}^N\|_{ E}^{q}\Big]&=
\E\Big[\big\|\sum_{i=0}^{k} T_{\delta t}^{k+1-i}G(Y^N_{k}+Z^N_k) \delta  W_i\big\|_{E}^{q}\Big]\\
&\le 
\E\Big[\Big(\sum_{i=0}^{k}\sum_{j\in \N^+}\big\| T_{\delta t}^{k+1-i}G(Y^N_{i}+Z^N_i)e_j\big\|_{\HH^{\frac 12+\epsilon}}^2\delta t\Big)^{\frac q2}\Big]\\
&\le \E\Big[\Big(\sum_{i=0}^{k}
((k+1-i)\delta t)^{-\frac 12-2\epsilon}(1+\|Y^N_{i}\|^{2\alpha}+\|Z^N_i\|^{2\alpha})\delta t\Big)^{\frac q2}\Big].
\end{align*}
The H\"older inequality yields that for any $l>2$, 
\begin{align*}
&\E[\|Z^N_{k+1}\|_{E}^q]\\
&\le
C\E\Big[(\sum_{i=0}^k((k+1-i)\delta t)^{-\frac {1+4\epsilon}2\frac {l}{(l-1)}}\delta t)^{(\frac {q(l-1)}{2{l}})}
(\sum_{i=0}^k(1+\|Y^N_i\|^{2\alpha l}+\|Z^N_i\|^{2\alpha l})\delta t)^{\frac {q}{2l}}\Big].
\end{align*}
Since $\alpha<1$, there always exists $l>2$ such that 
$2\alpha l<4$. Then  taking  such $l$ and 
combining the estimates \eqref{dis-pri-h-1}-\eqref{pri-znk0}, we 
obtain 
\begin{align}\label{dis-pri-z}
&\E[\|Z^N_{k+1}\|_{E}^q]\\\nonumber 
&\le
C(\alpha,T)
\E\Big[(\sum_{i=0}^k(1+\|Y^N_i\|^{4}+\|Z^N_i\|^{2\alpha l})\delta t)^{\frac {q}{2l}}\Big]\le C(X_0,T,q,\alpha).
\end{align}
\qed

With the help of the factorization method, we can get a stronger priori estimates of $Z^N$.

\begin{lm}\label{e-sup-znk}
Let $X_0\in \HH$, $T>0$ and $q\ge1$.
There exists a positive constant $ C(X_0,T,q)$
such that 
\begin{align*}
\E\Big[\sup_{k\le K}\|Z^N_k\|^{q}_{E}\Big]
&\le  C(X_0,T,q).
\end{align*}
\end{lm}

\textbf{Proof.}
We introduce a new auxiliary process $\widehat Z^N(t)$ defined by 
\begin{align*}
\widetilde Z^N(t)=\int_0^t S(t-s)G(Y^N_{[r]}+Z^N_{[r]})dW(r).
\end{align*}
Then by using the factorization method and the procedures in the proof of \cite{CH19}[Lemma 3.2] , it is not difficult to obtain that 
\begin{align*}
\E\Big[ \sup_{t\in[0,T]}\|\widetilde Z^N(t)\|_E^q\Big]&\le C(X_0,T,q).
\end{align*}
By the triangle inequality, we have that 
\begin{align}\label{tr-eq}
\E\Big[ \sup_{k\le K}\|Z^N_k\|_E^q\Big]
&\le C_q\E\Big[\sup_{k\le K}\|\widetilde Z^N(k\delta t)\|_E^q\Big]
+C_q\E\Big[\sup_{k\le K}\|\widetilde Z^N(k\delta t)-Z^N_k\|_E^q\Big].
\end{align}
Next, it suffices to show that 
\begin{align}\label{conv-z}
\sup_{k\le K}\E\Big[\|Z^N_k-\widetilde Z^N(k\delta t)\|^{q}_{E}\Big]\le C(X_0,T,q)\delta^{\eta q},
\end{align} 
for some $\eta>0$ and $q>1$ large enough. Indeed, from the proof of \cite[Corollary 4.9]{BCH18}, we can obtain that for $\eta q>1$,
\begin{align*}
\E\Big[\sup_{k\le K}\|\widetilde Z^N(k\delta t)-Z^N_k\|_E^q\Big]
&\le \sum_{k\le K}\E\Big[\|\widetilde Z^N(k\delta t)-Z^N_k\|_E^q\Big]\\
&\le KC(X_0,T,q)\delta t^{\eta q}\le C(X_0,T,q),
\end{align*}
which, together with \eqref{tr-eq}, yields the desired result.
Now, we prove \eqref{conv-z}.
By using the Sobolev embedding theorem, the Burkholder inequality,error estimate between $T_{\delta t}^k$ and $e^{-A^2k\delta t}$, it holds that for some $\epsilon>0$ and $\eta>0$,
\begin{align*}
&\E\Big[\|\widetilde Z^N(t_k)-Z_k^N\|_E^q\Big]\\
&\le C_q\E\Big[\|\widetilde Z^N(t_k)-Z_k^N\|_{\HH^{\frac 12+\epsilon}}^q\Big]\\
&\le C_q\E\Big[\Big\|\int_0^{t_k}(e^{-A^2(t_k-s)}-T_{\delta t}^{k-[s]})G(Z_{[s]}^N+Y_{[s]}^N)dW(r)\Big\|_{\HH^{\frac 12+\epsilon}}^q\Big]\\
&\le C_q \delta t^{\frac {\eta q}2}\E\Big[\Big(\int_0^{t_k}(1+(t_k-[s])^{-\frac 12-\frac \epsilon4-\eta})(1+\|Z_{[s]}^N\|^{2\alpha}+\|Y_{[s]}^N\|^{2\alpha})dr\Big)^{\frac q2}\Big].
\end{align*}
Applying the H\"older inequality and taking $\epsilon>0$ small enough, we get that for some  $l>2$  and small $\eta>0$ which depend on $\alpha$,
\begin{align*}
&\E\Big[\|\widetilde Z^N(t_k)-Z_k^N\|_E^q\Big]\\
&\le C_q \delta t^{\frac {\eta q}2}\Big(\int_0^{t_k}(1+(t_k-[s])^{-(\frac 12+\frac \epsilon 4+\eta)\frac {l}{l-1}})ds\Big)^{\frac {q(l-1)}{2l}}
\\
&\quad \times
 \E\Big[\Big(\int_0^{t_k}(1+\|Z_{[s]}^N\|^{2\alpha l}+\|Y_{[s]}^N\|^{2\alpha l})dr\Big)^{\frac q{2l}}\Big].
\end{align*}
The above estimate, together with \eqref{dis-pri-h-1} and \eqref{pri-znk0} yields that \eqref{conv-z}. Thus we complete the proof.
\qed

\begin{prop}\label{e-sup-xnk}
Let  $X_0\in \HH$, $T>0$, $q\ge1$.
There exists a positive constant $ C(X_0,T,q)$
such that 
\begin{align*}
\sup_{k\le K}\E\Big[\|X^N_k\|^{q}\Big]
&\le C(X_0,T,q).
\end{align*}
\end{prop}

\textbf{Proof.}
By Lemma \ref{e-sup-znk}, it suffices to bound $ Y_k^N$ in $\HH$.
Taking inner product on both sides of Eq. \eqref{ily} with $ Y^N_{k+1}$ in $\HH$ and using the Gagliardo--Nirenberg inequality,
we obtain 
\begin{align*}
\| Y^N_{k+1}\|^2
&\le \| Y^N_{k}\|^2
-(2-\epsilon)\delta t\|A Y^N_{k+1}\|^2\\
&\quad +C(\epsilon)
(1+\|Z^N_{k+1}\|_{E}^4+\|Y^N_{k+1}\|_{L^{4}}^{4})\delta t.
\end{align*}

The above estimate, together with \eqref{dis-pri-h-1}-\eqref{pri-znk0} and Proposition \ref{pri-znk}, yields that for $p\ge 1$,
\begin{align*}
&\| Y^N_{k+1}\|_{L^{2p}(\Omega; \mathbb H)}+\Big\|\sum_{k=0}^{K-1}\|A Y^N_{k+1}\|^2\delta t\Big\|_{L^p(\Omega;\mathbb R)}\\
&\le C(\epsilon, T)\Big\|\sum_{k=0}^{K-1}\|Z^N_{k+1}\|_{E}^4\delta t\Big\|_{L^p(\Omega;\mathbb R)} \\&\quad+C(\epsilon, T)\Big\|\sum_{k=0}^{K-1}\|Y^N_{k+1}\|_{L^{4}}^{4}\delta t\Big\|_{L^p(\Omega;\mathbb R)}+C(X_0,p,T)
\\
&\le C(X_0,p,T),
\end{align*}
which completes the proof.
\qed

\begin{rk}\label{rk-xn}
Similarly, one can obtain the a priori estimate in $E$.
More precisely, it holds that for $q\ge1$, $\gamma\le \frac 12$ and sufficient small $\epsilon>0$,
\begin{align*}
\E\big[\|X^N_k\|^q_{E}\big]
&\le C(X_0,T,q)(1+(k\delta t)^{-\frac q8+\frac {\gamma q}4-\epsilon q}),
\end{align*}
and for $\gamma > \frac 12$,
\begin{align*}
\E\big[\|X^N_k\|^q_{ E}\big]
&\le C(X_0,T,q).
\end{align*}
\end{rk}

In the following, we present the discrete continuity of 
the discrete convolution $\{Z^N_k\}_{k\le K}$ and the numerical approximation $\{X^N_k\}_{k\le K}$.

\begin{prop}\label{pri-znk}
Let $X_0\in \HH$, $T>0$, $q\ge1$. Then for all $\gamma\in (0,\frac 32)$, there exists $C(X_0,T,p,\gamma)>0$ such that  
\begin{align*}
\sup_{k\le K}\E\Big[\big\|Z^N_k\big\|_{\HH^{\gamma}}^q\Big]
&\le  C(X_0,T,p,\gamma)
\end{align*}
and 
\begin{align*}
\E\Big[\big\|Z^N_k-Z^N_{k_1}\big\|^q\Big]
&\le  C(X_0,T,p,\gamma)|(k-k_1)\delta t|^{\frac {\gamma q} 4},
\end{align*}
where $k_1,k\le K$.
\end{prop}

\textbf{Proof.}
By the Burkh\"older inequality and the smoothing effect of $T_{\delta t}$, we obtain that  
\begin{align*}
\E\Big[\|Z^N_k\|^{q}_{\mathbb H^{\gamma}}\Big]
&\le 
C_q\E\Big[\Big(\sum_{i=0}^{k}\sum_{j\in \N^+}\big\|T_{\delta t}^{k+1-i}G(Y^N_{k}+Z^N_k)e_j\big\|_{\mathbb H^{\gamma}}^2 \delta t\Big)^{\frac q2}\Big]\\
&\le 
 C_q\E\Big[\Big(\sum_{i=0}^{k}
((k+1-i)\delta t)^{-\frac 14-\frac \gamma 2}(1+\|Y^N_{i}\|^{2\alpha}+\|Z^N_i\|^{2\alpha})\delta t\Big)^{\frac q2}\Big].
\end{align*}
By similar arguments in the proof of Lemma \ref{sup-e-znk}, we complete the proof for the first estimate.

For convenience, we denote $[s]=\frac {[s]_{\delta t}} {\delta t}$ and $[s]_{\delta t}:=\max\{0,\delta t,\cdots,k\delta t,\cdots \}\cap [0,s]$.
Define the continuous extension of $Z^N_k$ as
\begin{align*}
\widehat Z^N(t)=\int_0^tT_{\delta t}^{-(\frac t{\delta t}-[s])}P^NG(X_{[s]}^N)dW(s),\; \widehat Z^N(0)=0,
\end{align*}
where $\widehat Z^N(t_k)=Z^N_k, t_k=k\delta t$.
Without loss of generality, we assume that $k>k_1$
The mild form of $\widehat Z^N$ and the Burkh\"older inequality yield that 
for $0\le s\le t\le T$, $q\ge 1$ and sufficient small $\epsilon_1>0$,
\begin{align*}
&\E\Big[\|\widehat Z^N(t_{k})-\widehat Z^N(t_{k_1})\|^q\Big]\\
&\le \E\Big[\Big\|\int_0^{t_{k_1}}\big(T_{\delta t}^{-(k-[r])}-T_{\delta t}^{-(k_1-[r])}\big)P^N\big(G(Y^N_{[r]}+\widehat Z^N([r]_{\delta t}))\big)dW(r)\Big\|^q\Big]\\
&\quad
+\E\Big[\Big\|\int_{t_{k_1}}^{t_{k}}T_{\delta t}^{-(k-[r])}P^NG(Y^N_{[r]}+\widehat Z^N([r]_{\delta t}))dW(r)\Big\|^q\Big]\\
&\le C\E\Big[\Big(\int_0^{t_{k_1}}\sum_{i\in \N^+}\Big\|T_{\delta t}^{-(k_1-[r])}\big(T_{\delta t}^{-(k-k_1)}-I\big)P^N\big(G(Y^N_{[r]}+\widehat Z^N([r]_{\delta t}))\big)e_i\Big\|^2dr\Big)^{\frac q2}\Big]\\
&\quad
+C\E\Big[\Big(\int_{t_{k_1}}^{t_{k}}\sum_{i\in \N^+}\Big\|T_{\delta t}^{-(k-[r])}P^NG(Y^N_{[r]}+\widehat Z^N([r]_{\delta t}))e_i\Big\|^2dr\Big)^{\frac q2}\Big]\\
&\le 
 C\E\Big[\Big(\int_0^{t_{k_1}} (k_1-[r]_{\delta t})^{-\frac \gamma 2-\frac 14-\epsilon_1}((k-k_1)\delta t)^{\frac \gamma 2}\|G(Y^N_{[r]}+\widehat Z^N([r]_{\delta t}))\|^2dr\Big)^{\frac q2}\Big]\\
&\quad
+C\E\Big[\Big(\int_{t_{k_1}}^{t_{k}}(k-[r]_{\delta t})^{-\frac 14-\epsilon}\|G(Y^N_{[r]}+\widehat Z^N([r]_{\delta t}))\|^2dr\Big)^{\frac q2}\Big].
\end{align*}
Then combining the above estimate and Proposition \ref{e-sup-xnk}, we complete the proof.
\qed

Next, we focus on the  regularity estimate of the numerical solution $\{X^N_k\}_{k\le K}$.

\begin{prop}\label{pri-ynk}
Let $X_0\in \HH^{\gamma}$, $\gamma \in (0,\frac 32)$, $T>0$ and  $q\ge1$.
Then there exists $ C(X_0,T,q,\gamma)>0$ such that
\begin{align*}
\sup_{k\le K}\E\Big[\big\|X^N_k\big\|_{\HH^{\gamma}}^q\Big]
&\le C(X_0,T,q,\gamma)
\end{align*}
and 
\begin{align*}
\E\Big[\big\|X^N_k-X^N_{k_1}\big\|^q\Big]
&\le C(X_0,T,q,\gamma)|(k-k_1)\delta t|^{\frac {\gamma p} 4},
\end{align*}
where $k_1,k\le K$.
\end{prop}

\textbf{Proof.}
Due to the above  regularity estimate in Lemma \ref{pri-znk}, 
it only suffices to give the regularity estimate of $\{Y^N_k\}_{k\le K}$.
To this end, we first give the upper bound of $\|Y_k^N\|_{L^6}$. For convenience, we only give the proof for the case that $\gamma\in [1,\frac 32)$.
From the Sobolev embedding theorem, the smoothing effect of $e^{-A^2t}$, the Gagliardo--Nirenberg and Young inequalities, it follows  that 
\begin{align*}
&\| Y^N_{k+1}\|_{L^6}\\
&\le \|T_{\delta t}^{k+1} Y_0^N\|_{L^6}
+ \|\sum_{j=0}^{k}T_{\delta t}^{k+1-j}A F(Y_{j+1}^N+Z_{j+1}^N)\|_{L^6}\delta t\\
&\le C\|Y_0^N\|_{\HH^1}+
C\sum_{j=0}^k(k+1-j)^{-\frac 7{12}}\delta t^{-\frac 7{12}}
\Big(1+\|Z_{j+1}^N\|_{L^{6}}^{3}+\| Y^N_{j+1}\|\Big)\delta t\\
&\quad +C\sum_{j=0}^k(k+1-j)^{-\frac 79}\delta t^{-\frac 79} \| Y^{N}_j\|^{\frac {10}3}\delta t+\sum_{j=0}^k\|(-A) Y^N_{j+1}\|^2\delta t.
\end{align*}
The a priori estimates of $Y^N_k$
 and $Z^N_k$ imply that for $q\ge 1$,
\begin{align*}
\E\Big[\sup_{k\le K}\| Y^N_{k}\|_{L^6}^q\Big]
&\le C(X_0,T,q).
\end{align*}
Now, we are in the position to give the desired regularity estimate. 
From the mild form of $Y^N_k$ and the above a priori estimate in $L^6$, it follows that 
\begin{align*}
&\E\Big[\big\| Y^N_k\big\|_{\HH^{\gamma}}^q\Big]\\
&\le 
C\E\Big[\| Y^N_0\|_{\HH^{\gamma}}^q\Big]
+C\delta t\E\Big[\Big(\sum_{j=0}^{K-1} \left\|T_{\delta t}^{K-j}A F(Y^N_{j+1}+Z^N_{j+1})\right\|_{\HH^{\gamma}}^q\Big]\\
&\le C(q)\| X^N(0)\|_{\HH^{\gamma}}^q+
C(q)\sum_{j=0}^{K-1}(T-t_j)^{-\frac 12-\frac \gamma 4}\delta t\sup_{j\le K}\E\Big[\Big(1+\|Y^N_{j}\|_{L^6}^3+\|Z^N_{j}\|_{L^6}^3\Big)^q\Big]\\
&\le C(q)\| X^N(0)\|_{\HH^{\gamma}}^q+
C(T,p)\sup_{j\le K}\E\Big[\Big(1+\|Y^N_{j}\|_{L^6}^3+\|Z^N_{j}\|_{L^6}^3\Big)^q\Big]\\
&\le C(X_0,T,q).
\end{align*}
For convenience, we assume that  $k> k_1$.
The mild forms of $Y^N_k$ and $Y^N_{k_1}$, together with  the smoothing effect of $T_{\delta t}$, yield that 
\begin{align*}
\E\Big[\big\|Y^N_k-Y^N_{k_1}\big\|^q\Big]
&\le \E\Big[\big\| T_{\delta t}^{k_1}(T_{\delta t}^{k-k_1}-I)Y_0^N\|^q\Big]
\\
&\quad+\E\Big[\Big(\sum_{j=0}^{k_1-1}\Big\|(T_{\delta t}^{k_1-j}(T_{\delta t}^{k-k_1}-I)A F(Y^N_{j+1}+Z^N_{j+1})\Big\|\delta t\Big)^{q}\Big]
\\
&\quad+\E\Big[\Big(\sum_{j=k_1}^{k-1}\Big\|T_{\delta t}^{k-j}AF(Y^N_{j+1}+Z^N_{j+1})\Big\|\delta t\Big)^q\Big]\\
&\le C(X_0,T,q,\gamma)|(k-k_1)\delta t|^{\frac {\gamma q} 4}.
\end{align*}
Combining the above regularity estimates together, we finish the proof for $\gamma\in [1,\frac 32)$.
Similar arguments, together with
 \begin{align*}
\|e^{-A^2t}X_0^N\|_{L^6}\le Ct^{-\frac 1{12}+\frac {\gamma }4}\|X_0\|,
\end{align*}
lead to the 
desired result for the case of $\gamma\in (0,1)$.
\qed

\subsection{Mean square convergence rate of the full discretization}
\label{sec-str}

Based on the a priori  and  regularity estimates of the numerical approximation, we focus on the strong convergence and the strong convergence rate of the proposed scheme for Eq. \eqref{spde}. 

For convenience, we only present the convergence analysis in mean square sense. 
Notice that $\|X^N_k-X(t_k)\|\le \|X^N_k-X^N(t_k)\|+\|X^N(t_k)-X(t_k)\|$. 
By Lemma \ref{lm-xn} and \cite[Remark 4.1]{CH19}, we have that 
if $ X_0\in \HH^{\gamma}$,  $\gamma \in (0,\frac 32)$, then for $\alpha\in (0,\gamma)$ it holds that 
\begin{align}\label{sup-err-h1}
&\|X^N-X\|_{L^2(\Omega;C([0,T];\HH))}
\le C(T,X_0,\alpha) \lambda_{N}^{-\frac \alpha 2}.
\end{align} 
for a positive constant $C(T,X_0,\alpha)$.
Thus we only need to estimate the error of $\|X^N_k-X^N(t_k)\|$.

To deduce the mean square convergence rate in time, 
we  introduce an auxiliary process $\widetilde X^N_k$, $k\le K$ with $\widetilde X^N_0=X^N(0)$, defined by 
\begin{align*}
\widetilde X^N_{k+1}
=\widetilde X^N_{k}
-A^2\delta t \widetilde X^N_{k+1}
-P^N A F(X^N(t_{k+1}))\delta t
+P^NG(X^N(t_{k}))\delta W_k.
\end{align*}
 Then we split the error of $X^N_{k}-X^N(t_{k})$ as 
\begin{align*}
\|X^N_{k}-X^N(t_{k})\| \le \|X^N(t_{k})-\widetilde X^N_{k}\| 
+\|\widetilde X^N_{k}-X^N_{k}\|.
\end{align*} 
The first error is bounded by the following lemma.
The second error will be dealt with the interpolation arguments. We would like to mention that compared with the existing literatures (see e.g. \cite{CH17,CHL16b,CHLZ17,HJ14}), this is also a new approach to deducing strong convergence rates of temporal and full discretizations for SPDEs with non-monotone coefficients.

\begin{lm}\label{lm-ynk}
Let $  X_0\in \HH^{\gamma}$, $\gamma \in (0,\frac 32)$ and  $p\ge1$.
There exists a positive constant $ C(X_0,T,p,\gamma)$ such that
\begin{align}\label{strong-ynk}
\big\|X^N(t_k)-\widetilde X^N_{k}\big\|_{L^p(\Omega;\HH)}
&\le C(X_0,T,p,\gamma)\delta t^{\frac \gamma 4}.
\end{align}
\end{lm}

\textbf{Proof.}
Denote $[s]_{\delta t}:=\max\{0,\delta t,\cdots,k\delta t,\cdots \}\cap [0,s]$ and $[s]=\frac {[s]_{\delta t}} {\delta t}$.
The mild forms of $X^N(t_{k})$ and $\widetilde X^N_{k}$ yield that 
\begin{align}\label{dec}
&\|X^N(t_{k})-\widetilde X^N_{k}\|_{L^p(\Omega;\HH)}\\\nonumber 
&\le
\Big\|\int_0^{t_{k}}(e^{-A^2(t_k-s)}-T_{\delta t}^{k-[s]})AP^NF(X^N(s))ds\Big\|_{L^p(\Omega;\HH)}\\\nonumber
&\quad+\Big\|\int_0^{t_{k}}T_{\delta t}^{k-[s]}AP^N\Big(F(X^N(s))-F(X^N({[s]_{\delta t}+\delta t}))\Big)ds\Big\|_{L^p(\Omega;\HH)}\\\nonumber
&\quad+\Big\|\int_0^{t_{k}}(e^{-A^2(t_k-s)}-T_{\delta t}^{k-[s]})P^NG(X^N(s))dW(s)\Big\|_{L^p(\Omega;\HH)}\\\nonumber
&\quad+\Big\|\int_0^{t_{k}}T_{\delta t}^{k-[s]}P^N\Big(G(X^N(s))-G(X^N({[s]_{\delta t}}))\Big)dW(s)\Big\|_{L^p(\Omega;\HH)} \\\nonumber
&=:I_1+I_2+I_3+I_4.
\end{align}
By the properties of $e^{-A^2t}$ and $T_{\delta t}^{k}$, the priori estimates of $Y^N$ and $Z^N$, and $\|e^{-A^2t}X_0^N\|_{L^6}\le Ct^{\min(-\frac 1{12}+\frac \gamma 4,0)}\|X_0\|$, it holds that for $\beta\in(0,\min(\frac 12,\frac 14+\frac {3\gamma}4))$,
\begin{align*}
I_1
&\le  C\sum_{j=0}^{k-1}\int_{t_j}^{t_{j+1}}(t_k-[s]_{\delta t})^{-\frac 12+\beta}\delta t^{\beta}
\Big(1+\Big\|X^N(s)\Big\|^3_{L^{3p}(\Omega;L^6)}\Big)ds\\
&\le C(T,X_0,p)\delta t^{\beta}.
\end{align*}
The similar arguments in the proof of \cite[Proposition 3.2]{CH19}  yield that 
for  $\epsilon>0$ small enough,
\begin{align*}
\E\big[\|X^N(t)\|^q_{E}\big]
&\le C(X_0,T,q)(1+t^{\min(0,-\frac q8+\frac {\gamma q}4-\epsilon q)}),
\end{align*}
which implies that 
\begin{align*}
I_{2}
&\le \Big\|\int_0^{t_{k}}T_{\delta t}^{k-[s]}AP^N\Big(F'(X^N(s))(X^N(s)-X^N({[s]_{\delta t}+\delta t}))\Big)ds\Big\|_{L^p(\Omega;\HH)}\\
&\le \int_0^{t_{k}} (t_k-[s]_{\delta t})^{-\frac 12} \Big\|\|F'(X^N(s))\|_{E}\|X^N(s)-X^N({[s]_{\delta t}+\delta t})\|\Big\|_{L^p(\Omega)}ds\\
&\le C(T,X_0)\delta t^{\frac \gamma 4}.
\end{align*}
The Burkh\"older inequality and the Sobolev embedding theorem lead that for small $\epsilon>0$,
{\small
\begin{align*}
I_3
&\le C(T,p) \Big(\Big\|\int_0^{t_{k}}\Big\|(e^{-A^2(t_k-s)}-T_{\delta t}^{k-[s]})P^NG(X^N(s))\Big\|_{\LL_2(\HH)}^2ds\Big\|_{L^{\frac p2}(\Omega)}\Big)^{\frac 12}\\
&\le C(T,p) \Big(\Big\|\int_0^{t_{k}}\big\|A^{\frac 14+\epsilon+\frac \gamma 2}(e^{-A^2(t_k-s)}-T_{\delta t}^{k-[s]})A^{-\frac \gamma 2}\big\|^2\\
&\qquad \qquad \Big\|A^{-\frac 14-\epsilon}P^NG(X^N(s))\Big\|_{\LL_2(\HH)}^2ds\Big\|_{L^{\frac p2}(\Omega)}\Big)^{\frac 12}\\
&\le C(T,p) \Big(\Big\|\int_0^{t_{k}}\big\|A^{\frac 14+\epsilon+\frac \gamma 2}(e^{-A^2(t_k-s)}-T_{\delta t}^{k-[s]})A^{-\frac \gamma 2}\big\|^2 \|A^{-\frac 14-\epsilon} G(X^N(s))\Big\|_{\LL_2(\HH)}^2ds \|_{L^{\frac p2}(\Omega)}\Big)^{\frac 12}\\
&\le C(X_0,T,p)\delta t^{\frac \gamma 4}. 
\end{align*}
}
Similar arguments yield that
{\small 
\begin{align*}
I_4
&\le C(T,p) \Big(\Big\|\int_0^{t_{k}}\Big\|T_{\delta t}^{k-[s]}P^N\big(G(X^N(s))-G(X^N({[s]_{\delta t}}))\big)\Big\|_{\LL_2(\HH)}^2 ds\Big\|_{L^{\frac p2}(\Omega)}\Big)^{\frac 12}\\
&\le  C(T,p)  \Big(\Big\|\int_0^{t_{k}}\big\|A^{\frac 14+\epsilon}T_{\delta t}^{k-[s]}\big\|^2 \Big\|A^{-\frac 14-\epsilon}\big(G(X^N(s))-G(X^N({[s]_{\delta t}}))\big)\Big\|_{\LL_2(\HH)}^2 ds\Big\|_{L^{\frac p2}(\Omega)}\Big)^{\frac 12}\\
&\le C(T,p)  \Big(\Big\|\int_0^{t_{k}} (t_k-[s]_{\delta t})^{-\frac 12-2\epsilon }\big\|X^N(s)-X^N({[s]_{\delta t}})\big\|^2 ds\Big\|_{L^{\frac p2}(\Omega)}\Big)^{\frac 12}\\
&\le C(X_0,T,p)  \delta t^{\frac \gamma 4}.
\end{align*}
}
Combining \eqref{dec} and the above regularity estimates, we complete the proof.
\qed

Next, we show the optimal regularity of  $\widetilde X^N_{k}$ and give
the estimate of $\widetilde X^N_{k}-X^N_k$ in $\HH^{-1}$. By a proof similar to that of \cite[Proposition 3.1]{CH19}, we have the following result on the regularity estimate of $\widetilde X^N_{k}$.

\begin{lm}\label{lm-pri-y}
Let $ X_0\in \HH^{\gamma}$, $\gamma \in (0,\frac 32)$, $T>0$ and $p\ge1$.
Then $\widetilde X^N_k$  satisfies 
\begin{align*}
\sup_{k\le K}\E\Big[\big\| \widetilde X^N_k\big\|_{\HH^{\gamma}}^p\Big]
&\le C(X_0,T,p,\gamma).
\end{align*}
\end{lm}

\begin{lm}\label{lm-ynk1}
Under the condition of Lemma \ref{lm-ynk}, there 
exist $\delta t_0\le 1$ and  $ C(X_0,T,\gamma)>0$ such that 
for any $\delta t\le \delta t_0$, $N\in \N^+$, we have 
\begin{align}\label{strong-ynk-h}
\big\|\widetilde X^N_{k}-X^N_k\big\|_{L^2(\Omega;\HH^{-1})}
&\le C(X_0,T,\gamma)\delta t^{\frac \gamma 4}.
\end{align} 
\end{lm}
\textbf{Proof.}
From the definitions of  $ X^N_{k+1}$ and $\widetilde X^N_{k+1}$, it follows that  for small $\epsilon>0$,
\begin{align*}
&\|X^N_{k+1}-\widetilde X^N_{k+1}\|_{\HH^{-1}}^2\\
&\le \|X^N_{k}-\widetilde X^N_{k}\|_{\HH^{-1}}^2-2\|A^{\frac 12}(X^N_{k+1}-\widetilde X^N_{k+1})\|^2 \delta t \\
&\quad
-2\<F(X^N_{k+1})
-F(X^N(t_{k+1})),X^N_{k+1}-\widetilde X^N_{k+1})\>\delta t
\\
&\quad+\<X^N_{k+1}-\widetilde X^N_{k+1}, (G(X_k^N)-G(X^N(t_k)))\delta W_k\>_{\HH^{-1}}\\
&\le \|X^N_{k}-\widetilde X^N_{k}\|_{\HH^{-1}}^2
-(2-\epsilon)\|A^{\frac 12}(X^N_{k+1}-\widetilde X^N_{k+1})\|^2 \delta t\\
&\quad+
2C(\epsilon)\|X^N_{k+1}-\widetilde X^N_{k+1}\|_{\HH^{-1}}^2\delta t
-2C(\epsilon)\|F(\widetilde X^N_{k+1})-F(X^N(t_{k+1}))\|^2_{\HH^{-1}}\delta t\\
&\quad+\<X^N_{k+1}-X^N_k, (G(X_k^N)-G(X^N(t_k)))\delta W_k\>_{\HH^{-1}}\\
&\quad+\<X^N_k-\widetilde X^N_{k}, (G(X_k^N)-G(X^N(t_k)))\delta W_k\>_{\HH^{-1}}
\\
&\quad+
\<\widetilde X_k^N-\widetilde X^N_{k+1}, (G(X_k^N)-G(X^N(t_k)))\delta W_k\>_{\HH^{-1}}\\
&=: \|X^N_{k}-\widetilde X^N_{k}\|_{\HH^{-1}}^2
-(2-\epsilon)\|A^{\frac 12}(X^N_{k+1}-\widetilde X^N_{k+1})\|^2 \delta t\\
&\quad+
2C(\epsilon)\|X^N_{k+1}-\widetilde X^N_{k+1}\|_{\HH^{-1}}^2\delta t
+II_1^k
+II_2^k+II_3^k+II_4^k.
\end{align*}
Taking $2C(\epsilon )\delta t$ small enough, for $l\le K-1$, we get
\begin{align}\label{iter}
&\|X^N_{l+1}-\widetilde X^N_{l+1}\|_{\HH^{-1}}^2+
\sum_{k=0}^{l}\|A^{\frac 12}(X^N_{k+1}-\widetilde X^N_{k+1})\|^2 \delta t
\\\nonumber
&\le 
2C(\epsilon)\sum_{k=0}^{l-1}\|X^N_{k+1}-\widetilde X^N_{k+1}\|_{\HH^{-1}}^2\delta t
+
C\sum_{k=0}^{l}(II_1^k
+II_2^k+II_3^k+II_4^k).
\end{align}
Taking expectation on both sides of \eqref{iter}, we have  
\begin{align*}
&\E\Big[\|X^N_{l+1}-\widetilde X^N_{l+1}\|_{\HH^{-1}}^{2}\Big]+
\E\Big[\sum_{k=0}^{l}\|A^{\frac 12}(X^N_{k+1}-\widetilde X^N_{k+1})\|^2 \delta t \Big]\\
&\le 
C(\epsilon)\sum_{k=0}^{l-1}\E\Big[\|X^N_{k+1}-\widetilde X^N_{k+1}\|_{\HH^{-1}}^{2}\Big]\delta t+
C\E\Big[\sum_{k=0}^{l}II_1^k\Big]\\
&\quad +C\E\Big[\sum_{k=0}^{l}II_2^k\Big]
+
C\E\Big[\sum_{k=0}^{l}II_3^k\Big]
+
C\E\Big[\sum_{k=0}^{l}II_4^k\Big].
\end{align*}

By using the a priori estimates of $X_{k+1}^N$, $X^N$, Lemma \ref{lm-ynk} and \cite[Lemma 3.5]{CH19}, we have that for small $\epsilon>0$,
\begin{align*}
&\E\Big[\sum_{k=0}^{l}II_1^k\Big]\\
&\le C\E\Big[\sum_{k=0}^{l}\| F(\widetilde X^N_{k+1})-F(X^N(t_{k+1}))\|_{\HH^{-1}}^2\delta t\Big]\\
&\le C(X_0,p)\E\Big[\sum_{j=0}^{l}(1+\|X^N(t_{j+1})\|_{E}^{4}+\|\widetilde X^N_{j+1}\|_{E}^{4})\|\widetilde X^N_{j+1}-X^N(t_{j+1})\|^{2}\delta t\Big]\\
&\le C(X_0,T,p)(\sum_{j=0}^l((j+1)\delta t)^{\min(-\frac 12+\gamma-\epsilon,0)}\delta t)\delta t^{\frac {\gamma } 2}.
\end{align*}
The mild form of $X^N_{k}$ yields that 
\begin{align*}
&\E\Big[\sum_{k=0}^{l}\Big(II_2^k+II_4^k\Big)\Big]\\
&\le
C\E\Big[\sum_{k=0}^l \<(T_{\delta t}-I)(X^N_k-\widetilde  X^N_k), (G(X_k^N)-G(X^N(t_k)))\delta W_k\>_{\HH^{-1}}\Big]
\\
&+C
\E\Big[\sum_{k=0}^l \<T_{\delta t}A(F(X^N_{k+1})-F(X^N(t_{k+1})))\delta t, (G(X_k^N)-G(X^N(t_k)))\delta W_k\>_{\HH^{-1}}\Big]\\
&+C
\E\Big[\sum_{k=0}^l \<T_{\delta t}(G(X^N_{k})-G(X^N(t_k)))\delta W_k, (G(X_k^N)-G(X^N(t_k)))\delta W_k\>_{\HH^{-1}}\Big]\\
&=:III_1+III_2+III_3.
\end{align*}
The Burkh\"older inequality and Sobolev interpolation inequality lead to 
\begin{align*}
III_1&\le 
\E\Big[\Big(\sum_{k=0}^l \sum_{j\in \N^+}\<(T_{\delta t}-I)(X^N_k-\widetilde X^N_k), (G(X_k^N)-G(X^N(t_k))e_j\>_{\HH^{-1}}^2\delta t\Big)^{\frac 12}\Big]\\
&\le 
\E\Big[\sup_{k\le l}\|(T_{\delta t}-I)(X^N_k-\widetilde X^N_k)\|_{\HH^{-1}} \Big(\sum_{k=0}^l \|X_k^N-\widetilde X_k^N\|^2\delta t\Big)^{\frac 12}
 \Big]\\
 &<\infty.
\end{align*}
Thus $III_1=0$.
For the second term, 
\begin{align*}
&III_2\\
&\le 
C\E\Big[\sum_{k=0}^l \<T_{\delta t}A(F(X^N_{k})-F(X^N(t_k)))\delta t, (G(X_k^N)-G(X^N(t_k)))\delta W_k\>_{\HH^{-1}}\Big]
\\
&+
C\E\Big[\sum_{k=0}^l \<T_{\delta t}A(F(X^N_{k+1})-F(X^N_k))\delta t, (G(X_k^N)-G(X^N(t_k)))\delta W_k\>_{\HH^{-1}}\Big]\\
&+
C\E\Big[\sum_{k=0}^l \<T_{\delta t}A(F(X^N(t_{k+1}))-F(X^N(t_k)))\delta t, (G(X_k^N)-G(X^N(t_k)))\delta W_k\>_{\HH^{-1}}\Big]\\
&=:III_{21}+III_{22}+III_{23}.
\end{align*}
Similar to $III_1$, we have $III_{21}=0$.
The estimations of  $III_{22}$ and $III_{23}$ are similar, we only give the estimate of $III_{22}$.
The continuity of $X^N_k$, Lemma \ref{lm-ynk} and  the Sobolev interpolation equality lead to
\begin{align*}
&III_{22}\\
&\le C\E\Big[\Big|\sum_{k=0}^l \<T_{\delta t}A(F(X^N_{k+1})-F(X^N_k))\delta t, (G(X_k^N)-G(X^N(t_k)))\delta W_k\>_{\HH^{-1}}\Big|\Big]
\\
&\le C\sum_{k=0}^l \E\Big[\|T_{\delta t}A(F(X^N_{k+1})-F(X^N_k))\delta t\|_{\HH^{-1}}^2\Big]\\
&\quad
+C\sum_{k=0}^l
\E\Big[\|(G(X_k^N)-G(X^N(t_k)))\delta W_k\|_{\HH^{-1}}^2\Big]\\
&\le C(T,X_0)\delta t^{\frac \gamma 2}+C\sum_{k=0}^l
\E\Big[\|X_k^N-\widetilde X^N_k\|^2\Big]\delta t
+C\sum_{k=0}^l
\E\Big[\|\widetilde X_k^N-X^N(t_k)\|^2\Big]\delta t\\
&\le C(T,X_0)\delta t^{\frac \gamma 2}+\epsilon \sum_{k=0}^l
\E\Big[\|X_k^N-\widetilde X^N_k\|_{\HH^1}^2\Big]\delta t
\\
&\quad+C(\epsilon) \sum_{k=0}^l
\E\Big[\|X_k^N-\widetilde X^N_k\|_{\HH^{-1}}^2\Big]\delta t
+C\sum_{k=0}^l
\E\Big[\|\widetilde X_k^N-X^N(t_k)\|^2\Big]\delta t.
\end{align*}
The Burkh\"older inequality  yields that 
\begin{align*} 
III_3&\le
C\sum_{k=0}^l \sum_{j\in \N^+}\E\Big[\|(G(X^N_{k})-G(X^N(t_k)))e_j\|_{\HH^{-1}}^2\Big]\delta t\\
&\le C(\epsilon)\sum_{k=0}^l\E\Big[\|X^N_{k}-X^N(t_k)\|_{\HH^{-1}}^2\Big] \delta t\\
&\quad+\epsilon \sum_{k=0}^l\E\Big[\|(X^N_{k}-X^N(t_k)\|_{\HH^1}^2\Big] \delta t.
\end{align*}
Notice that $\E\Big[\sum_{k=0}^{l}II_3^k\Big]=0$.
Combining all the estimates of $II_1^k$-$II_4^k$, we complete the proof.
\qed

\begin{prop}\label{tem-str}
Let $ X_0\in \HH^{\gamma}$,  $\gamma \in (0,\frac 32)$, $T>0$. Then there exists $\delta t_0\le 1$ such that for $N\in \N^+$ and $\delta t\in (0,\delta t_0)$,
\begin{align*}
\|X^N_k-X^N(t_k)\|_{L^2(\Omega; \HH)}
&\le C(T,X_0,\alpha)\delta t^{\frac \alpha 4},
\end{align*}
where  $\alpha\in (0,\gamma)$ and $C(T,X_0,p,\alpha)>0$. 
\end{prop}
\textbf{Proof.}
The triangle inequality yields that 
\begin{align*}
\|X^N_k-X^N(t_k)\|_{L^2(\Omega; \HH)}
&\le \|\widetilde X^N_k-X^N(t_k)\|_{L^2(\Omega; \HH)}
+\|X^N_k-\widetilde X^N_k\|_{L^2(\Omega; \HH)}\\
&\le C(T,X_0)\delta t^{\frac \gamma 4}
+\|X^N_k-\widetilde X^N_k\|_{L^2(\Omega; \HH)}.
\end{align*}
The mild forms of $X^N_l$ and $\widetilde X^N_l$
lead to
\begin{align*}
\|X^N_l-\widetilde X^N_l\|_{L^2(\Omega; \HH)}
&\le 
C\|\sum_{k=0}^{l-1}T_{\delta t}^{l-k}A(F(X^N_{k+1})-F(X^N(t_{k+1})))\delta t\|_{L^2(\Omega; \HH)}\\
&\quad+C\|\sum_{k=0}^{l-1}T_{\delta t}^{l-k}(G(X^N_k)-G(X^N(t_k)))\delta W_k\|_{L^2(\Omega; \HH)}. 
\end{align*}
The smoothing effect of $T_{\delta t}$, together with the Sobolev interpolation inequality, Remark \ref{rk-xn}, \cite[Lemma 3.5]{CH19} and Lemma \ref{lm-ynk1},
implies that for $\beta \in (0,1)$,
\begin{align*}
&\|\sum_{k=0}^{l-1}T_{\delta t}^{l-k}A(F(X^N_{k+1})-F(X^N(t_{k+1})))\delta t\|_{L^2(\Omega; \HH)}\\
&\le 
\sum_{k=0}^{l-1}(t_l-t_k)^{-\frac 34}\Big\|\Big(1+\|X^N_{k+1}\|_{ E}^2+\|X^N_{k+1}\|^2_{\HH^{\beta}}+\|X^N(t_{k+1})\|_{ E}^2\\
&\quad+\|X^N(t_{k+1})\|^2_{\HH^{\beta}}\Big) \|X^N_{k+1}-X^N(t_{k+1})\|^{\beta}_{\HH^{-1}}\|X^N_{k+1})-X^N(t_{k+1})\|^{1-\beta}\Big\|_{L^2(\Omega)}\delta t\\
&\le \sum_{k=0}^{l-1}(t_l-t_k)^{-\frac 34}
\Big\|\Big(1+\|X^N_{k+1}\|_{E}^2+\|X^N_{k+1}\|^2_{\HH^{\beta}}+\|X^N(t_{k+1})\|_{E}^2\\
&\quad  +\|X^N(t_{k+1})\|^2_{\HH^{\beta}}\Big)\|X^N_{k+1})-X^N(t_{k+1})\|^{1-\beta} \Big\|_{L^{\frac 2{1-\beta}}(\Omega)} \\
&\qquad \times \Big\|\|X^N_{k+1}-X^N(t_{k+1})\|^{\beta}_{\HH^{-1}}\Big\|_{L^{\frac 2\beta }(\Omega)}\delta t\\
&\le C(T,X_0)\delta t^{\frac {\beta \gamma} 4}.
\end{align*}
The Burkh\"older inequality yields that 
\begin{align*}
&\|\sum_{k=0}^{l-1}T_{\delta t}^{l-k}(G(X^N_k)-G(X^N(t_k)))\delta W_k\|^2_{L^2(\Omega; \HH)}\\
&\le \sum_{k=0}^{l-1}\sum_{j\in \N^+}\E[\|T_{\delta t}^{l-k}(G(X^N_k)-G(X^N(t_k)))e_j\|^2]\delta t\\
&\le \sum_{k=0}^{l-1} (t_k-t_l)^{-\frac 14}\E\Big[\|X^N_k-\widetilde X^N_k\|^2\Big] \delta t+C(T,X_0,\gamma)\delta t^{\frac \gamma 2}.
\end{align*}
The discrete Gronwall inequality leads to the desired result.
\qed

As a consequence of \eqref{sup-err-h1} and Proposition \ref{tem-str}, we present the sharp mean square convergence rate result of the proposed scheme. We also remark that the convergence analysis of the numerical approximation can be applied to more general $Q$-Wiener process cases. 

\begin{tm}\label{cor-main}
Let $ X_0\in \HH^{\gamma}$,  $\gamma \in (0,\frac 32)$, $T
>0$. Then there exists $\delta t_0\le 1$ such that for $N\in \N^+$ and $\delta t\in (0,\delta t_0)$,
\begin{align*}
\|X^N_k-X(t_k)\|_{L^2(\Omega; \HH)}
&\le C(T,X_0,\alpha)(\delta t^{\frac \alpha 4}+\lambda_N^{-\frac \alpha 2}),
\end{align*}
where $\alpha\in (0,\gamma)$ and $C(T,X_0,\alpha)>0$.
\end{tm}

From the above arguments in the proof of  Proposition \ref{tem-str}, it is not hard to see that under the same condition of Theorem \ref{cor-main}, $\Big\|\sup\limits_{k\le K}\|X^N_k-X(t_k)\|\Big\|_{L^2(\Omega,\R)}\le C(T,X_0,\alpha)(\delta t^{\frac \alpha 4}+\lambda_N^{-\frac \alpha 2})$, $\alpha \in (0,\gamma)$ also holds.
As another consequence of the above theorem, we have the following strong convergence of the proposed full discretization in 
$ E$-norm.
\begin{cor}\label{cor-co}
Let $ X_0\in \HH^{\gamma}$,  $\gamma \ge \frac 32$, $T
>0$. Then there exists $\delta t_0\le 1$ such that for $N\in \N^+$ and $\delta t\in (0,\delta t_0)$,
\begin{align*}
\|X^N_k-X(t_k)\|_{L^2(\Omega; E)}
&\le C(T,X_0,\alpha)(\delta t^{\frac \alpha 4}+\lambda_N^{-\frac \alpha 2}),
\end{align*}
where $\alpha\in (0,1)$ and $C(T,X_0,\alpha)>0$.
\end{cor}

The above result implies that once the density function of $X^N_k(x)$ exists, the density function of  $X^N_k(x)$ weakly converges or converges in Wasserstein metric
to  the density function of the exact solution $X(t_k,x)$. If one can further obtain the uniform Malliavin differentiability of $X^N_k(x)$, then the weak convergence of the density functions can be improved to the strong convergence.

\section{Conclusion}
In this article, we prove the 
the existence and positivity of the density function of the law of  the exact solution for the stochastic Cahn--Hilliard equation with unbounded diffusion coefficient satisfying a growth condition of order $0<\alpha<1$. Furthermore, we propose an implicit full discretization to approximate the density function of the exact solution. With the help of  the variation approach and the factorization method, the sharp mean square convergence rate of the full discretization is established. However, to fully understand how to compute the density function numerically, there are still many   open problems.
For instance,
 the existence and 
 positivity of the density function of full discretization, and the smoothness of the density functions of both the exact solution and full discretization are  still unknown for stochastic Cahn--Hilliard equation. 
\\

{\bf \large References}

\bibliographystyle{plain}
\bibliography{bib}

\end{document}